\documentclass[a4paper,10pt,leqno,twoside]{tobart}

\usepackage[english]{babel} \usepackage{inputenc, amsmath, amssymb, latexsym,
  epic, epsfig, rotating, fancyheadings, amsthm, pifont, empheq}

\newcommand{\follows}{\ensuremath{\Rightarrow}}

\newcommand{\ld}{\ensuremath{,\ldots,}}
\newcommand{\ssq}{\ensuremath{\subseteq}}
\newcommand{\smin}{\ensuremath{\setminus}}
\newcommand{\eps}{\ensuremath{\varepsilon}}

\newcommand{\T}{\ensuremath{\mathbb{T}}}

\newcommand{\N}{\ensuremath{\mathbb{N}}} 
\newcommand{\R}{\ensuremath{\mathbb{R}}}
\newcommand{\Z}{\ensuremath{\mathbb{Z}}}

\newcommand{\inte}{\ensuremath{\mathrm{int}}}

\newcommand{\diam}{\ensuremath{\mathrm{diam}}}

\newcommand{\kreis}{\ensuremath{\mathbb{T}^{1}}}

\newcommand{\torus}{\ensuremath{\mathbb{T}^2}}

\newcommand{\alphlist}{\begin{list}{(\alph{enumi})}{\usecounter{enumi}\setlength{\parsep}{2pt}
      \setlength{\itemsep}{1pt} \setlength{\topsep}{5pt}
      \setlength{\partopsep}{3pt}}}

\newcommand{\arablist}{\begin{list}{(\arabic{enumi})}{\usecounter{enumi}\setlength{\parsep}{2pt}
          \setlength{\itemsep}{1pt} \setlength{\topsep}{5pt}
          \setlength{\partopsep}{3pt}}}

\newcommand{\romanlist}{\begin{list}{(\roman{enumi})}{\usecounter{enumi}\setlength{\parsep}{2pt}
              \setlength{\itemsep}{1pt} \setlength{\topsep}{5pt}
              \setlength{\partopsep}{3pt}}}

 \newcommand{\listend}{\end{list}}

\newcommand{\bulletlist}{\begin{list}{$\bullet$}{\setlength{\parsep}{2pt}
                \setlength{\itemsep}{1pt} \setlength{\topsep}{5pt}
                \setlength{\partopsep}{3pt}\setlength{\leftmargin}{15pt}}}

\newcommand{\myproof}{\textit{Proof. }}

\newcommand{\foot}{\footnote}

\newcommand{\nLim}{\ensuremath{\lim_{n\rightarrow\infty}}}

\newcommand{\tLim}{\ensuremath{\lim_{t\rightarrow\infty}}}

\newcommand{\ntel}{\ensuremath{\frac{1}{n}}}

\newcommand{\halb}{\ensuremath{\frac{1}{2}}}

\newcommand{\viertel}{\ensuremath{\frac{1}{4}}}

\newcommand{\dreiviertel}{\ensuremath{\frac{3}{4}}}

\newtheoremstyle{tobthm}{3pt}{3pt}{\itshape}{0pt}{\bfseries}{.}{0.5eM}{}
\theoremstyle{tobthm}
\newtheorem{definition}{Definition}[section]
\newtheorem{thm}[definition]{Theorem}
\newtheorem{theorem}[definition]{Theorem}

\newtheorem{lemma}[definition]{Lemma}

\newtheoremstyle{tobrem}{3pt}{3pt}{\normalfont}{0pt}{\bfseries}{.}{0.5em}{}
\theoremstyle{tobrem} \newtheorem{remarks}[definition]{Remarks}
\newtheorem{rem}[definition]{Remark}

\newcommand{\eqand}{\ensuremath{\quad \textrm{and} \quad}}

\numberwithin{equation}{section}
\numberwithin{figure}{section}

\title{\Large\textsc{Nonautonomous saddle-node bifurcations: Random and
    deterministic forcing}} \author{V.~Anagnostopoulou \and T.~J\"ager}

\newcommand{\lmax}{\ensuremath{\lambda_{\mathrm{max}}}}

\begin{document}

\setlength{\abovedisplayskip}{0.8ex}
\setlength{\abovedisplayshortskip}{0.6ex}

\setlength{\belowdisplayskip}{0.8ex}
\setlength{\belowdisplayshortskip}{0.6ex}

\maketitle

\abstract{We  study the effect of external forcing on the saddle-node bifurcation
  pattern of interval maps. By replacing fixed points of unperturbed maps by
  invariant graphs, we obtain direct analogues to the classical result both for
  random forcing by measure-preserving dynamical systems and for deterministic
  forcing by homeomorphisms of compact metric spaces. Additional assumptions
  like ergodicity or minimality of the forcing process then yield further
  information about the dynamics.

  The main difference to the unforced situation is that at the critical
  bifurcation parameter, two alternatives exist. In addition to the possibility
  of a unique neutral invariant graph, corresponding to a neutral fixed point,
  a pair of so-called pinched invariant graphs may occur. In quasiperiodically
  forced systems, these are often referred to as `strange non-chaotic
  attractors'. The results on deterministic forcing can be considered as an
  extension of the work of Novo, N\'u\~nez, Obaya and Sanz on nonautonomous
  convex scalar differential equations. As a by-product, we also give a
  generalisation of a result by Sturman and Stark on the structure of minimal
  sets in forced systems. }

\noindent
\section{Introduction}

An important question which arises frequently in applications is that of the
influence of external forcing on the bifurcation patterns of deterministic
dynamical systems. This has been one of the main motivations for the development
of random dynamical systems theory (compare \cite[Chapter 9]{arnold:1998}), and the
description of the nonautonomous counterparts of the classical bifurcation
patterns is one of the principal goals of nonautonomous bifurcation theory. The
different types of forcing processes which are of interest range from deterministic systems
like quasiperiodic motion or, more generally, strictly ergodic dynamics on the one side
to random or stochastic processes like Brownian motion (white noise) at the
other end of the spectrum.  The reader is referred to \cite[Section 9]{arnold:1998} for a good
introduction to the topic and to
\cite{johnson/kloeden/pavani:2002,novo/obaya/sanz:2004,nunez/obaya:2007,homburg/young:2006,%
  zmarrou/homburg:2007,zmarrou/homburg:2008} for more recent developments and
further references.

Our aim here is to consider one of the simplest types of bifurcations, namely
saddle-node bifurcations of interval maps or scalar differential
equations. Given a forcing transformation $\omega : \Theta\to\Theta$, where
$\Theta$ is either a measure space or a topological space, we study skew product
maps of the form
\begin{equation} \label{e.skewproduct}
  f(\theta,x) \ : \ \Theta \times [a,b] \to \Theta \times [a,b] \quad ,
  \quad (\theta,x) \mapsto (\omega(\theta),f_\theta(x)) \ ,
\end{equation}
where $\omega:\Theta\to\Theta$ is called the {\em forcing process} or {\em base
  transformation}.  The bifurcating objects we concentrate on are {\em invariant
graphs}, that is, measurable functions $\varphi:\Theta \to [a,b]$ which satisfy
\begin{equation}
  f_\theta(\varphi(\theta)) \ = \ \varphi(\omega(\theta))
\end{equation}
for all (or at least almost all) $\theta\in\Theta$. Suppose we are given a
parameter family $(f_\beta)_{\beta\in[0,1]}$ of maps of the form
(\ref{e.skewproduct}) and a region $\Gamma\ssq\Theta\times[a,b]$. Then our
objective is to provide a criterium for the occurrence of saddle-node
bifurcations (of invariant graphs) inside of $\Gamma$. More precisely, we show
the existence of a critical bifurcation parameter $\beta_c$ such that
\begin{itemize}
\item If $\beta<\beta_c$, then $f_\beta$ has two invariant graphs in $\Gamma$.
\item If $\beta>\beta_c$, then $f_\beta$ has no invariant graphs in $\Gamma$.
\item If $\beta=\beta_c$, then $f_\beta$ has either one or two invariant graphs
  in $\Gamma$. If there exist two invariant graphs, then these are `interwoven'
  in a certain sense ({\em pinched}, Section~\ref{Pinched}).
\end{itemize}
Apart from some mild technical conditions, the crucial assumptions we need to
establish statements of this type are the monotonicity of the fibre maps
$f_\theta$, both with respect to $x$ and to the parameter $\beta$, and their
convexity inside of the considered region $\Gamma$ (see Theorems
\ref{t.random-sn} and \ref{t.deterministic-sn}).

Nonautonomous saddle-node bifurcations of this type have been studied previously in
\cite{novo/obaya/sanz:2004,nunez/obaya:2007} for nonautonomous scalar convex
differential equations over a strictly ergodic base flow and in
\cite{jaeger:2006a,nguyenetal:2010} for quasiperiodically forced interval
maps. In all cases, the proofs hinge on a convexity argument used to control the
number of invariant graphs or, more or less equivalently, minimal sets in the
system. This simple, but elegant and powerful idea can be traced back to Keller
\cite{keller:1996} and has later been used independently by Alonso and Obaya
\cite{alonso/obaya:2003} in order to classify nonautonomous scalar convex
differential equations according to the structure of their minimal
sets. However, so far no systematic use of these arguments has been made in
order to determine the greatest generality to which the description of
nonautonomous saddle-node bifurcations can be pushed. This is the goal of the
present paper.
Quite surprisingly, it turns out that hardly any assumptions on the underlying
forcing process are needed in order to give a fairly good description of the
bifurcation pattern. We only require that the forcing
transformation is invertible and that it is either a measure-preserving
transformation of a probability space or a homeomorphism of a compact metric
space. In the former case, we work in a purely measure-theoretic setting, such
that no topological structure on the base space is required. Additional
properties like ergodicity, respectively minimality, can be used in order
to obtain further information about the dynamics.

As a by-product of our studies in the topological setting, we also obtain a
generalisation of a result by Sturman and Stark \cite{sturman/stark:2000}
concerning the structure of invariant sets. If a compact invariant set of a
minimally driven ${\cal C}^1$-map on a Riemannian manifold only admits negative
upper Lyapunov exponents (with respect to any invariant measure supported on
$M$), then $M$ is just a finite union of continuous curves (see
Theorem~\ref{Generalization of Theorem 14}).  \medskip

The paper is organised as follows. In Section~\ref{Preliminaries}, we collect a
number of preliminaries on forced interval maps, including the convexity result
due to Keller. In Section~\ref{Pinched}, we introduce and discuss various
concepts of inseparability of invariant graphs ({\em pinching}), which are
variations of the well-known notion of pinched sets and graphs for
quasiperiodically forced monotone interval maps
\cite{glendinning:2002,stark:2003}. Section~\ref{RandomForcing} then contains
the bifurcation result for randomly forced systems. In
Section~\ref{SturmanStark}, we provide the above-mentioned generalisation of
Sturman and Stark's result and use it in Section~\ref{DeterministicForcing} to
prove the bifurcation result for deterministic forcing. In
Section~\ref{ContinuousTime}, we
discuss the application to continuous-time systems and the relations to the
respective results of \cite{novo/obaya/sanz:2004,nunez/obaya:2007}. Finally, in 
Section~\ref{Examples}, we present some explicit examples to illustrate the results.
\medskip

\noindent{\bf Acknowledgements.} This work was supported by the
Emmy-Noether-grant `Low-dimensional and Nonautonomous Dynamics' (Ja
1721/2-1) of the German Research Council.

\section{Invariant measures, invariant graphs and Lyapunov
  exponents} \label{Preliminaries}

Given a transformation $\omega:\Theta\to\Theta$ of a {\em base space} $\Theta$, an {\em
  $\omega$-forced map} is a skew-product map 
\begin{equation}
  \label{e.w-forced} f:\Theta \times X \to \Theta \times X \quad , \quad  (\theta,x)
  \ \mapsto \ (\omega(\theta),f_\theta(x)) \ .
\end{equation}
$X$ is called the {\em phase space} and the maps $f_\theta :X \to X$ are called
{\em fibre maps}.  By $f^n_\theta=(f^n)_\theta$ we denote the fibre maps of the
iterates of $f$ (and not the iterates of the fibre maps). We will mostly
consider two situations: First, we study the case where $\Theta$ is a measurable
space, equipped with a $\sigma$-algebra ${\cal B}$, and $\omega$ is a measurable
bijection that leaves invariant a probability measure $\mu$.\foot{In all of the
  following, `measure' refers to a probability measure, unless explicitly
  stated otherwise.} This means that $(\Theta,{\cal B},\mu,\omega)$ is a
measure-preserving dynamical system, in the sense of Arnold \cite{arnold:1998},
with time $\T=\Z$. Secondly, we will treat the case where $\Theta$ is a compact
metric space and $\omega$ is a homeomorphism. In this case we always equip
$\Theta$ with the Borel $\sigma$-algebra ${\cal B}(\Theta)$. Consequently, for
any $\omega$-invariant Borel measure $\nu$ we arrive at situation one by taking
${\cal B}={\cal B}(\Theta)$ and $\mu=\nu$. However, it is important to emphasise
that we will not {\em a priori} fix any particular invariant measure in this
second setting. $X$ will always be a Riemannian manifold and in most cases
simply a compact interval $X=[a,b]\ssq \R$.

In the context of forced systems, fixed points of unperturbed maps are replaced
by {\em invariant graphs}. If $\mu$ is an $\omega$-invariant measure and $f$ is
an $\omega$-forced map, then we call a measurable function $\varphi:\Theta\to X$
an {\em $(f,\mu)$-invariant graph} if it satisfies
\begin{equation}
  \label{e.inv-graph}
  f_\theta(\varphi(\theta)) \ = \ \varphi(\omega(\theta)) \quad \textrm{for }
  \mu\textrm{-a.e. } \theta\in\Theta \ .
\end{equation}
When (\ref{e.inv-graph}) holds for all $\theta\in\Theta$, we say $\varphi$ is an
{\em $f$-invariant graph}, and in this case it is certainly an
$(f,\mu)$-invariant graph for all $\omega$-invariant measures $\mu$.  Usually,
we will only require that $(f,\mu)$-invariant graphs are defined $\mu$-almost
surely, which means that implicitly we always speak of equivalence
classes. Conversely, $f$-invariant graphs are always assumed to be defined
everywhere. This is particularly important in the topological setting, since in
this case topological properties like continuity or semi-continuity of the
invariant graphs play a role, and these can easily be destroyed by modifications
on a set of measure zero. As an additional advantage, the definition becomes
independent of an invariant reference measure on the base, which may not be
unique in the topological setting as we have mentioned before.

We say $f$ is an {\em $\omega$-forced monotone ${\cal C}^r$-interval map} if
$X=[a,b]\ssq \R$ and all fibre maps $f_\theta$ are $r$ times continuously
differentiable and strictly monotonically increasing. When $\omega$ is a
continuous map, we assume in addition that all derivatives $f_\theta^{(k)}(x),\
k=0\ld r$, depend continuously on $(\theta,x)$. The {\em (vertical) Lyapunov
  exponent} of an $(f,\mu)$-invariant graph $\varphi$ is given by
\begin{equation}
  \label{e.lyap}
  \lambda_\mu(\varphi) \ = \ \int_\Theta \log f'_\theta(\varphi(\theta)) \ d\mu(\theta) \ .
\end{equation}
For $\omega$-forced monotone interval maps with convex fibre maps, the following
result allows to control the number of invariant graphs and their Lyapunov
exponents at the same time.

\begin{thm}[Keller \cite{keller:1996}] \label{t.convexity} Let $(\Theta,{\cal
    B},\mu,\omega)$ be a {\em mpds} and $f$ be an $\omega$-forced ${\cal
    C}^2$-interval map. Further, assume there exist measurable functions
  $\gamma^-\leq \gamma^+ : \Theta\to X$ such that for $\mu$-a.e.\
  $\theta\in\Theta$ the maps $f_\theta$ are strictly monotonically increasing
  and strictly convex on
  $\Gamma_\theta=[\gamma^-(\theta),\gamma^+(\theta)]$. Further, assume that the
  function $\eta(\theta) = \inf_{x\in I(\theta)} \log f'_\theta(x)$ has an
  integrable minorant.

  Then there exist at most two $(f,\mu)$-invariant graphs in
  $\Gamma=\left\{(\theta,x)\in\Theta\times X \mid \gamma^-(\theta) \leq x \leq
    \gamma^+(\theta)\right\}$.\foot{We say an $(f,\mu)$-invariant graph
    $\varphi$ is contained in $\Gamma$ if there holds
    $\varphi(\theta)\in\Gamma_\theta$ $\mu$-a.s.~.} Further, if there exist two
  distinct $(f,\mu)$-invariant graphs $\varphi^-\leq \varphi^+$ in $\Gamma$ then
  $\lambda_\mu(\varphi^-)<0$ and $\lambda_\mu(\varphi^+)>0$.
\end{thm}
Implicitly, this result is contained in \cite{keller:1996}. A proof in the
quasiperiodically forced case, which literally remains true in the more general
situation stated here, is given in \cite{jaeger:2006a}.

Apart from the analogy to fixed points of unperturbed maps, an important reason
for concentrating on invariant graphs is the fact that there is a one-to-one
correspondence between invariant graphs and invariant ergodic measures of forced
monotone interval maps. On the one hand, if $f$ is an $\omega$-forced map, $\mu$
is an $\omega$-invariant ergodic measure and $\varphi$ is an $(f,\mu)$-invariant
graph, then an $f$-invariant ergodic measure $\mu_\varphi$ can be defined by
\begin{equation}
  \label{e.graph-measure}
  \mu_\varphi(A) \ = \ \mu\left(\left\{\theta\in\Theta\mid (\theta,\varphi(\theta)) \in A
   \right\}\right) \ .
\end{equation}
Conversely, we have the following.
\begin{thm}[Theorem 1.8.4 in \cite{arnold:1998}] \label{t.measures}
  Suppose $(\Theta,{\cal B},\mu,\omega)$ is an ergodic {\em mpds} and
  $f$ is an $\omega$-forced monotone ${\cal C}^0$-interval map.
  Further, assume that $\nu$ is an $f$-invariant ergodic measure which
  projects to $\mu$ in the first coordinate. Then $\nu=\mu_\varphi$
  for some $(f,\mu)$-invariant graph $\varphi$.
\end{thm}
The proof in \cite{arnold:1998} is given for the continuous-time case, but the
adaption to the discrete-time setting is immediate.

\section{Pinched invariant graphs} \label{Pinched}

An important notion in the context of minimally forced one-dimensional maps is
that of pinched sets and pinched invariant graphs
\cite{glendinning:2002,stark:2003,jaeger/stark:2006,fabbri/jaeger/johnson/keller:2005}. In
order to introduce it, we need some more notation. Let $X=[a,b]\ssq \R$.  Given
two measurable functions $\varphi^-,\varphi^+ :\Theta\to X$, we let
$$[\varphi^-,\varphi^+]=\left\{(\theta,x)\mid
  x\in\left[\varphi^-(\theta),\varphi^+(\theta)\right]\right\} \ , $$ similarly
for open and half-open intervals. For a subset $A\ssq \Theta\times X$ with
$\pi_1(A)=\Theta$, we let
\begin{equation}\label{e.bounding-graphs}
\varphi^-_A(\theta) = \inf A_\theta \quad \textrm{and} \quad
\varphi^+_A(\theta)=\sup A_\theta \ ,
\end{equation}
where $A_\theta = \{ x\in X\mid (\theta,x)\in A\}$.  Note that when $\Theta$ is
a topological space and $A$ is compact, then $\varphi^+_A$ is lower
semi-continuous (l.s.c.) and $\varphi^-_A$ is upper semi-continuous (u.s.c.).
Given $\varphi:\Theta\to X$, we denote the point set
$\Phi:=\{(\theta,\varphi(\theta))\mid \theta\in\Theta\}$ by the corresponding
capital letter. We let $\varphi^\pm := \varphi^\pm_{\overline{\Phi}}$ and write
$\varphi^{+-}$ and $\varphi^{-+}$ instead of $(\varphi^+)^- =
\varphi^-_{\overline{\Phi}^+}$ and $(\varphi^-)^+ =
\varphi^+_{\overline{\Phi}^-}$, ect.~.
\begin{definition}[Pinched graphs]\label{d.pinched}
  Suppose $\Theta$ is a compact metric space, $X=[a,b]\ssq\R$,
  $\varphi^-:\Theta\to X$ is l.s.c., $\varphi^+:\Theta\to X$ is u.s.c. and
  $\varphi^-\leq \varphi^+$. Then $\varphi^-$ and $\varphi^+$ are called {\em
    pinched} if there exists a point $\theta\in\Theta$ with
  $\varphi^-(\theta)=\varphi^+(\theta)$.

  A compact subset $A\ssq \Theta\times X$ with $\pi_1(A)=\Theta$ is called
  {\em pinched} if $\varphi^-_A$ and $\varphi^+_A$ are pinched, that is, if there
  exists some $\theta\in\Theta$ with $\# A_\theta=1$.
\end{definition}
There is a close relation between pinched graphs and minimal sets.
\begin{lemma}[\cite{stark:2003}]\label{l.stark}
  Suppose $\omega$ is a minimal homeomorphism of a compact metric space and $f$
  is an $\omega$-forced monotone ${\cal C}^0$-interval map. Then the following
  hold.  \alphlist
\item If $\varphi^-$ and $\varphi^+$ are pinched semi-continuous $f$-invariant
  graphs, then there exists a residual set $R\ssq \Theta$ with
  $\varphi^-(\theta)=\varphi^+(\theta) \ \forall \theta\in R$.
\item Any $f$-minimal set $A$ is pinched.
\item Any pinched compact $f$-invariant set $A$ contains exactly one minimal
  set.  \listend
\end{lemma}
The proof in \cite{stark:2003} is given for the case of quasiperiodic forcing,
but literally goes through for minimally forced maps. A slightly weaker concept
of pinching is the following.
\begin{definition}[Weakly pinched graphs] \label{d.weakly_pinched} Suppose
  $\Theta$ is a compact metric space, $X=[a,b]\ssq\R$, $\varphi^-:\Theta\to X$
  is l.s.c., $\varphi^+:\Theta\to X$ is u.s.c. and $\varphi^-\leq
  \varphi^+$. Then $\varphi^-$ and $\varphi^+$ are called {\em weakly pinched}
  if $\inf_{\theta\in\Theta} \varphi^+(\theta)-\varphi^-(\theta)=0$. Otherwise,
  we call $\varphi^-$ and $\varphi^+$ {\em uniformly separated}.
\end{definition}
Note that when $\varphi^-$ and $\varphi^+$ are uniformly separated, then there
exists some $\delta>0$ with $\varphi^-(\theta)\leq \varphi^+(\theta)-\delta \
\forall \theta\in\Theta$.

In the case of random forcing, a measure-theoretic analogue of pinching is
required.
\begin{definition}[Measurably pinched graphs]
  Suppose $(\Theta,{\cal B},\mu)$ is a measure space, $X = [a,b]\ssq \R$ and
  $\varphi^-\leq\varphi^+ : \Theta \to X$ are measurable. Then $\varphi^-$ and
  $\varphi^+$ are called {\em measurably pinched}, if the set
  $A_\delta:=\left\{\theta\in\Theta \mid \varphi^+(\theta)-\varphi^-(\theta) <
    \delta\right\}$ has positive measure for all $\delta>0$.  Otherwise, we call $\varphi^-$ and
  $\varphi^+$ {\em $\mu$-uniformly separated}.
  \end{definition}
  Similar to above, when $\varphi^-$ and $\varphi^+$ are $\mu$-uniformly
  separated there exists $\delta>0$ with $\varphi^-(\theta) \leq
  \varphi^+(\theta)-\delta$ for $\mu$-a.e.\ $\theta\in\Theta$.  In the case of
  minimal forcing, all three notions of pinching coincide.
  \begin{lemma} \label{e.pinching-equivalence}
    Suppose $\omega$ is a minimal homeomorphism of a compact metric space
    $\Theta$ and $f$ is an $\omega$-forced monotone ${\cal C}^0$-interval map.
    Further, assume that $\varphi^-\leq\varphi^+:\Theta\to X$ are $f$-invariant
    graphs, with $\varphi^-$ l.s.c and $\varphi^+$ u.s.c.~.

    Then $\varphi^-$ and $\varphi^+$ are pinched if and only they are weakly
    pinched if and only they are measurably pinched with respect to every
    $\omega$-invariant measure $\mu$ on $\Theta$.
 \end{lemma}
 \myproof We first show that pinching implies measurable pinching. Suppose that
 $\varphi^-$ and $\varphi^+$ are pinched, $\mu$ is an $\omega$-invariant measure
 and $\delta >0$. Then the set $A_\delta=\left\{\theta\in\Theta \mid
   \varphi^+(\theta)-\varphi^-(\theta) < \delta\right\}$ is non-empty and open
 (openness follows from the semi-continuity of $\varphi^\pm$). By minimality
 $\Theta=\bigcup_{i=0}^k \omega^{-i}(U)$ for some $k\in\N$. Then, by the
 $\omega$-invariance of $\mu$, $\mu(A_\delta)>0$. As $\delta>0$ was arbitrary,
 $\varphi^-$ and $\varphi^+$ are measurably pinched.

 The fact that measurable pinching implies weak pinching is obvious. Hence, in
 order to close the circle, assume that $\varphi^-$ and $\varphi^+$ are weakly
 pinched. Suppose for a contradiction that $\varphi^-$ and $\varphi^+$ are not
 pinched, such that $P=\{\theta \in\Theta \mid
 \varphi^-(\theta)=\varphi^+(\theta)\}$ is empty. Let $A_n=
 \{\theta\in\Theta\mid \varphi^+(\theta)-\varphi^-(\theta) \geq 1/n\}$. As
 $\Theta\smin P = \bigcup_{n\in\N} A_n$ is a countable union of closed sets,
 Baire's Theorem implies that for some $n\in\N$ the set $A_n$ has non-empty
 interior. Let $U=\inte(A_n)$. By minimality $\Theta=\bigcup_{i=0}^k
 \omega^i(U)$ for some $k\in\N$. The uniform continuity of $f$ on $\Theta\times
 X$ implies that there exists some $\delta>0$, such that $|x-y|\geq 1/n$ implies
 $|f^i_\theta(x)-f^i_\theta(y)| \geq \delta$ for all $\theta\in\Theta$ and
 $i=0\ld k$. Due to the invariance of the graphs $\varphi^\pm$ we therefore obtain
 $\varphi^+(\theta)-\varphi^-(\theta) \geq \delta \ \forall \theta\in\Theta$, in
 contradiction to the definition of weak pinching.  \qed \medskip

\section{Saddle node bifurcations: Random forcing} \label{RandomForcing}

In this section we suppose that $(\Theta,{\cal B},\mu,\omega)$ is a {\em mpds}
and consider parameter families $(f_\beta)_{\beta\in[0,1]}$ of $\omega$-forced
monotone ${\cal C}^2$-interval maps
$f_\beta(\theta,x)=(\omega(\theta),f_{\beta,\theta}(x))$. In order to show that
these families undergo a saddle-node bifurcation, we need to impose a number of
conditions. These will be formulated in a semi-local way, meaning that we do not
make assumptions on the whole space $\Theta\times X$. Instead, we restrict our
attention to a subset $\Gamma=[\gamma^-,\gamma^+]$, with measurable functions
$\gamma^-\leq \gamma^+:\Theta\to X$, and describe bifurcations of invariant
graphs contained in $\Gamma$. Consequently, all the required conditions only
concern the restrictions of the fibre maps $f_\theta$ to the intervals
$\Gamma_\theta = \left[\gamma^-(\theta),\gamma^+(\theta)\right]$.  One advantage
of this formulation is that it allows to describe local bifurcations taking
place in forced non-invertible interval maps. We shall not pursue this issue
further here, but refer the interested reader to \cite{nguyenetal:2010}, where
this idea is used to describe the creation of 3-periodic invariant graphs in the
quasiperiodically forced logistic map.

\begin{theorem}[Saddle-node bifurcations, random forcing] \label{t.random-sn}
  Let $(\Theta,{\cal B},\mu,\omega)$ be a measure-preserving dynamical system and
  suppose that $(f_\beta)_{\beta\in [0,1]}$ is a parameter family of
  $\omega$-forced ${\cal C}^2$-interval maps. Further, assume that there exist
  measurable functions $\gamma^-,\gamma^+ : \Theta \to X$ with
  $\gamma^-<\gamma^+$ such that the following hold (for $\mu$-a.~e.\
  $\theta\in\Theta$ and all $\beta\in[0,1]$ where applicable). \romanlist
\item[(r1)] There exist two $\mu$-uniformly separated $(f_0,\mu)$-invariant
  graphs, but no $(f_1,\mu)$-invariant graph in $\Gamma$;
 \item[(r2)] $f_{\beta,\theta}(\gamma^\pm(\theta)) \geq
   \gamma^\pm(\omega(\theta))$;
 \item[(r3)] the maps $(\beta,x)\mapsto f_{\beta,\theta}(x)$ and
 $(\beta,x)\mapsto f'_{\beta,\theta}(x)$ are continuous;
 \item[(r4)] the function $\eta(\theta) = \sup\left\{|\log f'_{\beta,\theta}(x)| \mid
   x\in \Gamma_\theta,\ \beta\in[0,1]\right\}$ is integrable with respect to $\mu$;
 \item[(r5)] $f'_{\beta,\theta}(x) > 0 \ \forall x\in \Gamma_\theta$;
 \item[(r6)] there exist constants $0< c_1\leq C$ such that
   $c_1\leq \partial_\beta f_{\beta,\theta}(x) \leq C \ \forall x \in
   \Gamma_\theta$;
 \item[(r7)] there exists a constant $c_2>0$ such that $f''_{\beta,\theta}(x) >
   c_2 \ \forall x\in \Gamma_\theta$.
  \listend Then there exist a unique
   critical parameter $\beta_\mu\in(0,1)$ such that: \begin{itemize}
   \item If $\beta<\beta_\mu$ then there exist exactly two $(f_\beta,\mu)$-invariant
     graphs $\varphi^-_{\beta}<\varphi^+_{\beta}$ in $\Gamma$ which are $\mu$-uniformly
     separated and satisfy $\lambda(\varphi^-_{\beta})<0$ and $\lambda(\varphi^+_{\beta})>0$.
   \item If $\beta=\beta_\mu$ then either there exists exactly one
     $(f_\beta,\mu)$-invariant graph $\varphi_{\beta}$ in $\Gamma$, or
     there exist two $(f_\beta,\mu)$-invariant graphs
     $\varphi^-_{\beta}\leq \varphi^+_{\beta}$ in $\Gamma$ which are
     measurably pinched. In the first case
     $\lambda_\mu(\varphi_{\beta})=0$, in the second case
     $\lambda_\mu(\varphi^-_{\beta})< 0$ and
     $\lambda_\mu(\varphi^+_{\beta})>0$.
\item If $\beta > \beta_\mu$ then there are no $(f_\beta,\mu)$-invariant graphs
  in $\Gamma$.
\end{itemize}
\end{theorem}
\begin{rem} \label{r.random}
  It may seem surprising at first sight that there always exists a
  unique bifurcation parameter in the above situation, despite the
  possible lack of ergodicity. However, this uniqueness is due to the
  fact that we require invariant graphs to be defined over the whole
  base space. Taking into account invariant graphs which are only
  defined over $\omega$-invariant subsets of $\Theta$ yields a whole
  spectrum of bifurcation parameters, one for each $\omega$-invariant
  subset, and in this sense uniqueness does require ergodicity. We
  discuss these issues in detail after the proof of
  Theorem~\ref{t.random-sn}.
\end{rem}
\begin{remarks} \label{r.random-forcing}
  \alphlist
\item We denote the critical bifurcation parameter by $\beta_\mu$ in order to
  keep the dependence on $\mu$ explicit. This will become important in the
  topological setting of Section~\ref{DeterministicForcing}, where we do not a
  priori fix a particular invariant reference measure, but have to take
  different measures into account.
\item Assumptions (r1)--(r4) should be considered as rather mild technical
  conditions. The crucial ingredients are the monotonicity in $x$ (r5), the
  monotonicity in $\beta$ (r6) and the convexity of the fibre maps (r7).
\item The generality concerning the forcing process is surely optimal, with the
  only exception of infinite measure preserving processes which are not
  considered here. In particular, $\omega$ may simply be taken the identity.  In
  this case the fibre maps become independent monotone interval maps, and
  $\beta_\mu$ is the last parameter for which a saddle-node bifurcation has only
  occurred for a set of $\theta$'s of measure zero.

  In contrast to this, we leave open the question whether the strong uniform
  assumptions concerning the behaviour on the fibres can be weakened under
  additional assumptions on the forcing process, for example when the forcing is
  ergodic.
\item Symmetric versions of the above result hold for parameter families with
  concave fibre maps and/or with decreasing behaviour on the parameter
  $\beta$. These versions can be derived from the above one by considering the
  coordinate change $(\theta,x)\mapsto(\theta,-x)$ and the parametrisation
  $\beta\mapsto 1-\beta$.
\item The information on the Lyapunov exponents allows to describe the behaviour
  of almost-all points for $\beta\leq\beta_\mu$: For $\mu$-a.e.\
  $\theta\in\Theta$ all points between $\varphi^-_\beta(\theta)$ and
  $\varphi^+_\beta(\theta)$ converge to the lower graph, in the sense that
  $\lim_{n\to\infty}
  \left|f^n_{\beta,\theta}(x)-\varphi^-_\beta(\omega^n(\theta))\right|=0$.
  Points below $\varphi^-$ converge to $\varphi^-$ in the same sense, whereas
  all points above $\varphi^+$ eventually leave $\Gamma$ (compare \cite[Proposition
  3.3 and Corollary 3.4]{jaeger:2003}).  \listend
\end{remarks}

\proof[Proof of Theorem~\ref{t.random-sn}] We start with some preliminary
remarks and fix some notation. First, note that we may assume without
loss of generality that the fibre maps $f_{\beta,\theta}$ are strictly monotonically
increasing on all of $X$ and thus invertible. Otherwise $f_\beta$ can be modified
outside $\Gamma$ accordingly. This does not change the dynamics in $\Gamma$
and therefore does not affect the number and properties of the invariant graphs
contained in this set.

Given an $\omega$-forced monotone interval map $f$ and a measurable function
$\gamma$, we define its {\em forwards} and {\em backwards graph transforms}
$f_*\gamma$ and $f^{-1}_*\gamma$ by
\begin{equation} \label{e.graph-transform} f_*\gamma(\theta) \ := \
  f_{\omega^{-1}(\theta)}(\gamma(\omega^{-1}(\theta))) \quad \textrm{and} \quad
  f^{-1}_*\gamma(\theta) \ := \  f^{-1}_{\omega(\theta)}(\gamma(\omega(\theta))) \ .
\end{equation}
Further, we define sequences
\begin{equation} \label{e.gamma_n}
\gamma^-_{\beta,n}\ := \ f^n_{\beta
  *}\gamma^- \quad \textrm{and} \quad \gamma^+_{\beta,n}\ :=\ f^{-n}_{\beta*}\gamma^+ \ .
\end{equation}
Due to (r2) and (r5) the sequence $\gamma^-_{\beta,n}$ is increasing and
$\gamma^+_{\beta,n}$ is decreasing. Obviously, if there exists an
$(f,\mu)$-invariant graph in $\Gamma$ then both sequences remain bounded in
$\Gamma$ and thus converge pointwise to limits
\begin{equation}
  \label{e.gamma-limits}
  \varphi^-_\beta \ := \ \nLim \gamma^-_{\beta,n} \eqand
   \varphi^+_\beta \ := \ \nLim \gamma^+_{\beta,n} \ .
\end{equation}
Using the continuity of the fibre maps $f_{\beta,\theta}$ it is easy to see that
$\varphi^\pm_\beta$ are $(f_\beta,\mu)$-invariant graphs. More precisely,
$\varphi^+_\beta$ is the highest and $\varphi^-_\beta$ is the lowest
$(f_\beta,\mu)$-invariant graph in $\Gamma$.

In fact, in order to ensure the existence of invariant graphs in $\Gamma$ it
suffices to have a measurable function $\psi :\Theta\to X$ with
$\psi(\theta)\in \Gamma_\theta \ \forall \theta\in\Theta$ and $f_{\beta*}\psi
\leq \psi$. In this case the sequence $\gamma_{\beta,n}^-$ remains bounded in
$\Gamma$ since $\gamma^-\leq \gamma_{\beta,n}^- \leq f_{\beta*}^n\psi \leq \psi
\leq \gamma^+\ \forall n\in\N$, such that again $\varphi^-_\beta$ in
(\ref{e.gamma-limits}) (and consequently also $\varphi^+_\beta$) defines an
invariant graph. In particular, in this situation
\begin{equation}
  \label{e.intermediate-graph}
  \varphi^-_\beta \ \leq \ f_*\psi \ \leq \ \psi \ \leq \ \varphi^+_\beta \ .
\end{equation}

We now define the critical parameter by
  \begin{equation}\beta_\mu = \sup\left\{\beta\in[0,1] \left| \
 \forall \beta'<\beta \ \exists \textrm{ 2
          uniformly separated } (f,\mu)\textrm{-invariant graphs}\right.\right\} \ .
\end{equation}
\underline{$\beta<\beta_\mu$:} \quad By definition, there exist two uniformly
separated $(f,\mu)$-invariant graphs for all $\beta <
\beta_\mu$. Theorem~\ref{t.convexity} implies that these are the only ones and
that their Lyapunov exponents have the right signs. \medskip

\noindent \underline{$\beta>\beta_\mu$:} \quad Suppose that $\beta>\beta_\mu$
and there exists an $(f_\beta,\mu)$-invariant graph $\psi$ in $\Gamma$. Then
(r6) implies that for any $\beta'<\beta$ we have
  \begin{equation}
    \label{e.r4}
    f_{\beta'*}\psi \ \leq \ \psi -\eta\ ,
  \end{equation}
  where $\eta:= (\beta-\beta')\cdot c_1$. Hence, (\ref{e.intermediate-graph})
  implies that
  \begin{equation}
    \varphi^-_{\beta'} \ \leq \ f_*\psi \ \leq
    \ \psi - \eta \ \leq \ \varphi^+_{\beta'}-\eta \ .
  \end{equation}
  Consequently $f_{\beta'}$ has two uniformly separated $(f,\mu)$-invariant
  graphs for all $\beta' < \beta$, contradicting the definition of
  $\beta_\mu$.\medskip

\noindent\underline{$\beta=\beta_\mu$:} \quad
By the above reasoning, the two uniformly separated $(f_\beta,\mu)$-invariant
graphs for $\beta<\beta_\mu$ are $\varphi^\pm_\beta$ defined in
(\ref{e.gamma-limits}). Due to (r6), $\varphi^-_\beta$ increases as $\beta$ is
increased, whereas $\varphi^+_\beta$ decreases (since this is true for the
sequences $\gamma^-_{\beta,n}$ and $\gamma^+_{\beta,n}$, respectively). In
particular, as $\beta\nearrow\beta_\mu$ the two sequences converge $\mu$-almost
surely to graphs $\tilde\varphi^-$ and $\tilde\varphi^+$. These graphs
are $(f_{\beta_\mu},\mu)$-invariant, since
\begin{eqnarray*}
  \lefteqn{\left|f_{\beta_\mu,\theta}(\tilde\varphi^\pm(\theta)) -
      \tilde\varphi^\pm(\omega(\theta)) \right| \  \leq }
  \\ & &  \underbrace{\left|f_{\beta_\mu,\theta}(\tilde\varphi^\pm(\theta)) -
      f_{\beta,\theta}(\varphi_{\beta}^\pm(\theta))\right|}_{\stackrel{n\to\infty}{\longrightarrow}\
    0 \textrm{ by (r3)}}
  \ + \    \underbrace{\left| \varphi^\pm_\beta(\omega(\theta))-
      \tilde\varphi^\pm(\omega(\theta))\right|}_{\stackrel{n\to\infty}{\longrightarrow}\
    0 \textrm{ by definition of } \varphi^\pm_{\beta_\mu}}
  \ \longrightarrow \ 0 \quad (\textrm{as } \beta \nearrow \beta_\mu).
\end{eqnarray*}
We have $\tilde\varphi^\pm = \lim_{\beta\to\beta_\mu} \varphi^\pm_\beta =
\lim_{\beta\to\beta_\mu} \lim_{n\to\infty} \gamma^\pm_{\beta,n}$, and due to the
monotonicity of the sequences we may exchange the two limits on the right to
obtain $\tilde\varphi^\pm = \varphi^\pm_{\beta_\mu}$.

We claim that either either $\varphi^-_{\beta_\mu}=\varphi^+_{\beta_\mu}$
$\mu$-a.s. or $\varphi^-_{\beta_\mu}$ and $\varphi^+_{\beta_\mu}$ are measurably
pinched.  The only alternative to this is that $\varphi^-_{\beta_\mu}$ and
$\varphi^+_{\beta_\mu}$ are $\mu$-uniformly separated. In this case let
$\psi(\theta)
=(\varphi^+_{\beta_\mu}(\theta)-\varphi^-_{\beta_\mu}(\theta))/2$. We now use
the following elementary lemma.
\begin{lemma}
  Suppose $g:X\to X$ is ${\cal C}^2$ with $g'>0$ and $g''>c_2$ and let
  $\delta>0$. Then there exists a constant $\eps=\eps(c_2,\delta)$ such that for
  all $x,y\in X$ with $d(x,y)\geq \delta$ there holds
  \begin{equation}
    g\left(\frac{x+y}{2}\right) \ \leq \ \frac{g(x)+g(y)}{2} - \eps \ .
  \end{equation}
\end{lemma}
Since $\varphi^-_{\beta_\mu}$ and $\varphi^+_{\beta_\mu}$ are $\mu$-uniformly
separated and the fibre maps $f_{\beta_\mu,\theta}$ are uniformly convex by
(r7), it follows that for some $\eps > 0$ there holds $f_{\beta_\mu *}\psi \leq
\psi-\eps$. This together with (r6) implies that for all $\beta\leq
\beta_\mu+\frac{\eps}{2C}$ there holds $f_{\beta*}\psi \leq \psi
-\frac{\eps}{2}$. From (\ref{e.intermediate-graph}) we now obtain that
\begin{equation}
  \varphi^-_\beta \ \leq \ f_{\beta*}\psi \ \leq \ \psi -\frac{\eps}{2} \ \leq
  \ \varphi^+_\beta \quad \quad \forall \beta \in\left[\beta_\mu,\beta_\mu+\frac{\eps}{2C}\right] \ .
\end{equation}
Hence for all $\beta\in\left[\beta_\mu,\beta_\mu+\frac{\eps}{2C}\right]$ the
graphs $\varphi^-_\beta$ and $\varphi^+_\beta$ are $\mu$-uniformly separated, in
contradiction to the definition of $\beta_\mu$.

It remains to prove the statement about the Lyapunov exponents. When
$\varphi^-_{\beta_\mu}$ and $\varphi^+_{\beta_\mu}$ do not belong to the same
equivalence class, then $\lambda_\mu(\varphi^-_{\beta_\mu}) < 0$ and
$\lambda_\mu(\varphi^+_{\beta_\mu})>0$ follow from Theorem~\ref{t.convexity}.
Further, we have
\begin{eqnarray*}
  \lambda_\mu(\varphi^\pm_{\beta_\mu}) & = & \int_\Theta \log
  f_{\beta_\mu,\theta}'(\varphi^\pm_{\beta_\mu}(\theta)) \ d\mu(\theta)\\
  &  = & \ \lim_{\beta\nearrow\beta_\mu}  \int_\Theta \log
  f_{\beta,\theta}'(\varphi^\pm_{\beta}(\theta)) \ d\mu(\theta)\ =
   \ \lim_{\beta\nearrow\beta_\mu} \lambda_\mu(\varphi^\pm_\beta) \ .
  \end{eqnarray*}
  For the second equality, note that
  \[
  \log f_{\beta,\theta}'(\varphi^\pm_{\beta}(\theta)) \
  \ \stackrel{\beta\nearrow\beta_\mu}{\longrightarrow}\  \log
  f_{\beta_\mu,\theta}'(\varphi^\pm_{\beta_\mu}(\theta))
  \] pointwise due to (r3), and by (r4) we can apply dominated convergence with
  majorant $\eta$.

  This implies that $\lambda_\mu(\varphi^-_{\beta_\mu}) \leq 0$ and
  $\lambda_\mu(\varphi^+_{\beta_\mu}) \geq 0$, and when both graphs are
  $\mu$-a.s.\ equal their common Lyapunov exponent must therefore be zero.
\qed

\medskip\medskip We close this section with some remarks on the restriction of
the dynamics to invariant subsets, which mostly concerns the case of non-ergodic
forcing. Suppose $M$ is an $\omega$-invariant subset of $\Theta$ of positive
measure. Let $\mu_M(A)=\mu(A\cap M)/\mu(M)$ be the induced probability measure
on $M$. Then Theorem \ref{t.random-sn} holds for the measure-preserving
dynamical system $(M,{\cal B},\mu_M,\omega_{|M})$ and the parameter family
$f_{\beta|M\times X}$ with new bifurcation parameter
\begin{equation*}
  \beta_\mu^M = \sup\left\{\beta\in[0,1] \left| \ \forall \beta'<\beta \
      \exists \textrm{ 2 uniformly separated } (f_{\beta}\mid_{M\times
        X},\mu_M)\textrm{-invariant graphs}\right.\right\}.
\end{equation*}
Obviously, we have
\begin{rem} \label{r.parameter-monotonicity} Let $M\subset\Theta$ be
  such that $\omega(M)=M$ and $\mu(M)\in(0,1]$. Then
  $\beta^M_\mu\geq\beta_\mu$.
\end{rem}
Consequently, invariant graphs defined on subsets of $\Theta$ may
still exist after the bifurcation parameter $\beta_\mu$. For
simplicity of exposition, it is convenient to extend the
definition in (\ref{e.gamma-limits}) in the following way.
\begin{displaymath}
\begin{array}{cc}
  \varphi_\beta^-(\theta)=\left\{\begin{array}{ll}
      \displaystyle\nLim\gamma_{\beta,n}^-(\theta) &
  \textrm{, if }\gamma_{\beta,n}^-(\theta)\in\Gamma_\theta\forall n\\
      +\infty & \textrm{, otherwise}
                  \end{array}\right. ,
                & \varphi_\beta^+(\theta)=\left\{\begin{array}{ll}
                    \displaystyle\nLim\gamma_{\beta,n}^+(\theta) &
 \textrm{, if }\gamma_{\beta,n}^+(\theta)\in\Gamma_\theta\forall n\\
                    -\infty & \textrm{, otherwise}
                  \end{array}\right. .
\end{array}
\end{displaymath}
By (r6) $\beta\mapsto\gamma_{\beta,n}^-(x)$ is increasing for all
$x\in\Gamma_{\theta}$.  Further, it is easy to check that (r6) implies that
$\beta\mapsto f^{-1}_{\beta,\theta}(x)$ is decreasing, and hence
$\beta\mapsto\gamma_{\beta,n}^+(x)$ is decreasing for all $x\in\Gamma_{\theta}$.
This yields the following lemma.
\begin{lemma} \label{l.graph-monotonicity}
  For $\mu$-almost all $\theta\in\Theta$ the function
  $\beta\mapsto\varphi_\beta^-(x)$ is increasing and the function
  $\beta\mapsto\varphi_\beta^+(x)$ is decreasing.
\end{lemma}
We call an orbit ${\cal O}(\theta,x)$ $\Gamma$-bounded if $f_\beta^n(\theta,x)\in\Gamma
\ \forall n\in\Z$. The next lemma highlights the connection between invariant
graphs and $\Gamma$-bounded orbits.
\begin{lemma} \label{l.bounded-orbits}
  Consider the set of $\Gamma$-bounded orbits $$K(\beta) \ := \ \{ (\theta,x)\in\Theta
  \times X \mid {\cal O}(\theta,x) \textrm{ is } \Gamma\textrm{-bounded} \}$$ and its projection
  $B(\beta):= \pi_1(K(\beta))$. Then the following hold for all $\beta\in[0,1]$.
\romanlist\item
$K(\beta)$ is $f_\beta$-invariant, $B(\beta)$ is $\omega$-invariant.
\item
$K(\beta)=[\varphi_\beta^-,\varphi_\beta^+]$.
\item
If $\beta>\beta'$, then $K(\beta)\subseteq K(\beta')$ and $B(\beta)\subseteq B(\beta')$.
\listend
\end{lemma}
\proof (i) is obvious. For (ii), note that since
$[\varphi_\beta^-,\varphi_\beta^+]$ is $f_\beta$-invariant it follows
that $[\varphi_\beta^-,\varphi_\beta^+]\subseteq K(\beta)$.  Now let
$(\theta,x)\in\Gamma\backslash[\varphi_\beta^-,\varphi_\beta^+]$ and
assume first that $x>\varphi_\beta^+(\theta)$.  Then
$x>\gamma_{\beta,n}^+(\theta)$ for some $n\in\N$,
i.e.~$x>f^{-n}_{\beta,\omega^n(\theta)}(\gamma^+(\omega^n(\theta)))$.
Using (r5) we see that
$f_{\beta,\theta}^n(x)>\gamma^+(\omega^n(\theta))$, such that
$f_{\beta}^n(\theta,x)\notin\Gamma$ and therefore $(\theta,x) \notin
K(\beta)$. The case where $x<\varphi_\beta^+(\theta)$ is treated
similarly.

Now (iii) follows from (ii) since the invariant graphs
$\varphi_\beta^-$, $\varphi_\beta^+$ are increasing, respectively decreasing
with $\beta$ by Lemma~\ref{l.graph-monotonicity}.  \qed\medskip

In light of the preceeding statement, we can define a second `last'
bifurcation parameter
\[
\hat\beta_\mu \ := \ \sup\{\beta\in[0,1] \mid \mu(B(\beta))>0\} \
\]
and a bifurcation interval $I_\mu=[\beta_\mu,\hat\beta_\mu]$ over which the set
of $\Gamma$-bounded orbits vanishes. The case where $\omega$ is the identity easily
allows to produce examples where this happens in a continuous way over a
non-trivial interval. Note also that $\mu(B(\hat \beta_\mu))$ may or may not be
zero.

If $\omega$ is ergodic, then the fact that $B(\beta)$ is $\omega$-invariant
implies that $K(\beta)$ vanishes immediately.
\begin{lemma}
If $\omega$ is ergodic, then $\mu(B(\beta))=1$ for
$\beta\leq\beta_\mu$, and $\mu(B(\beta))=0$ for $\beta>\beta_\mu$.
\end{lemma}

\section{The existence of continuous invariant graphs} \label{SturmanStark}

The purpose of this section is to provide criteria, in terms of Lyapunov
exponents, which ensure that a compact invariant set $K$ of a forced ${\cal
  C}^1$-map consists of a finite union of continuous
curves. Lemma~\ref{l.curve-criterium} below treats the relatively simple case of
driven interval maps. This statement is crucial for passing from the
measure-theoretic setting in Section~\ref{RandomForcing} to the topological one
in Section~\ref{DeterministicForcing} below and will be a key ingredient in
the proof of Theorem~\ref{t.deterministic-sn}. Because of its intrinsic
interest, we also include a generalisation that holds for forced
${\cal C}^1$-maps on Riemannian manifolds, provided that the forcing
homeomorphism is minimal (Theorem~\ref{Generalization of Theorem 14} below).
This extends a result for quasiperiodically forced systems by Sturman and Stark
\cite{sturman/stark:2000}.

\begin{lemma} \label{l.curve-criterium} Suppose $\omega$ is a homeomorphism of a
  compact metric space $\Theta$, $f$ is an $\omega$-forced ${\cal C}^1$-interval
  map and $K$ is a compact $f$-invariant set that intersects every fibre
  $\{\theta\}\times X$ in a single interval, that is,
  $K=[\varphi^-_K,\varphi^+_K]$. Further, assume that for all $\omega$-invariant
  measures and all $(f,\mu)$-invariant graphs $\varphi$ contained in $K$ we have
  $\lambda_\mu(\varphi)< 0$. Then $K$ is just a continuous $f$-invariant curve.
\end{lemma}
For the proof, we need the following semi-uniform ergodic theorem from
\cite{sturman/stark:2000}.  Given a measure-preserving transformation $T$ of a
probability space $(Y,{\cal B},\nu)$ and a subadditive sequence of integrable
functions $g_n : Y \to \R$ (that is, $g_{n+m}(y) \leq g_n(y) + g_m(T^ny)$), the
limit
\[
\bar g(y) \ = \ \nLim g_n(y)/n
\]
exists $\nu$-a.s.\ by the Subadditive Ergodic Theorem
(e.g.~\cite{arnold:1998,katok/hasselblatt:1997}). Furthermore $\bar g$ is
$T$-invariant. Consequently, when $T$ is ergodic then $\bar g$ is
$\nu$-a.s.\ equal to the constant $\nu(\bar g) = \int_Y \bar g \
d\nu$.

\begin{theorem}[Theorem 1.12 in \cite{sturman/stark:2000}]\label{Semi-uniform ergodic}
  Suppose that $T:Y\rightarrow Y$ is a continuous map on a compact metrizable
  space $Y$ and $g_n : Y \to \R \ (n\in\N_0)$ is a subadditive sequence of
  continuous functions. Let $\tau$ be a constant such that $\nu(\bar g)<\tau$ for
  every $T$-invariant ergodic measure $\nu$. Then there exist $\delta>0$ and
  $N\in\N$, such that
\[
\frac{1}{N} \sum_{n=0}^{N-1}g_n(y)\leq \tau-\delta \qquad
\forall y\in Y\ .
\]
\end{theorem}

\proof[Proof of Lemma~\ref{l.curve-criterium}] Due to Theorem~\ref{t.measures},
any $f$-invariant ergodic measure $\nu$ is of the form $\nu=\mu_\varphi$ for
some $\omega$-invariant ergodic measure $\mu$ and an $(f,\mu)$-invariant graph
$\varphi$. Consequently, we have
\begin{equation}
  \label{e.curve-lemma1}
  \int_{\Theta\times X} \log f'_\theta(x) \ d\nu(\theta,x) \ = \
  \int_\Theta \log f'_\theta(\varphi(\theta)) \ d\mu(\theta) \ = \ \lambda_\mu(\varphi) \ < \ 0 \ .
\end{equation}
Hence, Theorem~\ref{Semi-uniform ergodic} with $Y=\Theta\times X$,
$T=f$, $\tau=0$ and $g_n(\theta,x) = \log (f^n_\theta)'(x)$ implies
that for some $N\in\N$ and $\alpha\in(0,1)$ we have
\begin{equation}
  \left(f^N_\theta\right)'(x) \ \leq \ \alpha \quad \forall (\theta,x) \in K \ .
\end{equation}
If we let $C := \left(\sup_{\theta\in\Theta}
  \varphi^+(\theta)-\varphi^-(\theta)\right)$, then this implies
\begin{equation}
  \diam\left(K_\theta\right) \ = \ \diam\left(f^N_{\omega^{-N}\left(\theta\right)}
  \left(K_{\omega^{-1}\left(\theta\right)}\right)\right) \ \leq \ \alpha \cdot
  \diam\left(K_{\omega^{-N}\left(\theta\right)}\right) \ \leq \ \alpha \cdot C
  \quad \forall \theta\in\Theta \ ,
\end{equation}
which yields $C\leq \alpha\cdot C$. This means that $C=0$, such that $K$ is the graph
of the continuous function $\varphi^-\equiv\varphi^+$.
\qed\medskip

When the underlying homeomorphism $\omega$ is minimal, then a similar statement
holds in much greater generality, namely for arbitrary compact invariant sets of
$\omega$-forced ${\cal C}^1$-maps on any Riemannian manifold. For the case of
quasiperiodic forcing by an irrational rotation of the circle, this was shown by
Sturman and Stark \cite[Theorem~1.14]{sturman/stark:2000}.  Their proof should
generalise to irrational rotations on higher-dimensional tori, but in any case
it makes strong use of the fact that the forcing transformation $\omega$ is an
isometry and of the existence of a smooth structure on $\Theta$. In contrast to
this, we want to consider the general case of a minimal base transformation
$\omega$ on an arbitrary compact metric space $\Theta$. The argument we present
below allows to bypass the technical problems due to weaker hypotheses on
$\Theta$ and also significantly reduces the length the proof. 

In the remainder of this section we let $X$ be a Riemannian manifold, endowed
with the canonical distance function $d$ induced by the Riemannian metric. We
suppose $f$ is an $\omega$-forced ${\cal C}^1$-map on $\Theta\times X$. The
upper Lyapunov exponent of $(\theta,x)\in\Theta\times X$ is
\begin{equation}
  \label{e.pointwise-lyap}
  \lmax(\theta,x) \ = \ \limsup_{n\to\infty} \ntel \log \| Df_\theta^n(x) \| \ ,
\end{equation}
where $Df_\theta(x)$ is the derivative matrix of $f_\theta$ in $x$ and $\| \cdot \|$
denotes the usual matrix norm.  Given any $f$-invariant probability measure
$\nu$, we define the upper Lyapunov exponent of $\nu$ by
\begin{equation}
  \label{e.maximal-Lyap}
  \lambda_{\textrm{max}}(\nu) \ = \ \int \lmax(\theta,x) \ d\nu(\theta,x) \ .
\end{equation}
Further, we let $X_k = \{x\in X^k \mid x_i \neq x_j \textrm{ if } i\neq j\}$ and
endow $X_k$ with the Hausdorff topology.
\begin{thm}\label{Generalization of Theorem 14}
  Suppose $\omega:\Theta\to\Theta$ is a minimal homeomorphism, $X$ is a
  Riemannian manifold, $f$ is an $\omega$-forced ${\cal C}^1$-map on
  $\Theta\times X$ and $K$ is a compact invariant set of $f$. Further, assume
  that $\lmax(\nu)<0$ for all $f$-invariant ergodic measures $\nu$ supported on
  $K$. Then there exist $k\in\N$ and a continuous map $\psi:\Theta \mapsto X_k$
  such that $K$ is the graph of $\psi$, that is,
$$
K= \big\{(\theta,\psi_i(\theta))\big|\ \theta\in\Theta,\ i=1\ld k \big\}\ .
$$
\end{thm}
\begin{rem}
  \alphlist
\item Note that since we do not assume any specific structure on $\Theta$, it
  does not make sense to speak of the smoothness of the curve $\psi$ in this
  setting (in contrast to \cite{sturman/stark:2000}). However, when $\Theta$ is
  a torus and $\omega$ and irrational rotation, then the smoothness of $\psi$ follows
  from its continuity \cite{stark:1999}. In general, smoothness can only be
  expected when $\omega$ is an isometry.
\item If $f$ is invertible, as in the case of forced monotone interval maps, the
  conclusion of Theorem \ref{Generalization of Theorem 14} also holds if
  $\lmax(\nu)>0$ for all ergodic measures $\nu$.  \listend
\end{rem}
\proof Applying Theorem~\ref{Semi-uniform ergodic} to $Y=\Theta\times X$, $T=f$,
$\tau=0$ and $\varphi_n(\theta,x) = \log\|Df^n_\theta(x)\|$, we obtain that for
some $N\in\N$ and $\alpha'\in(0,1)$
\begin{equation}
  \label{e.3-1}
  \|Df^N_\theta(x)\| \ \leq \ \alpha' \quad \forall (\theta,x) \in K \ .
\end{equation}
Replacing $f$ by $f^N$ if necessary, we may assume without loss of generality
$N=1$. By compacity, there exist some $\eps > 0$ and $\alpha \in (\alpha',1)$
such that
\begin{equation}
  \label{e.3-2}
  \|Df_\theta(x)\| \ \leq \ \alpha \quad \forall (\theta,x) \in B_\eps(K) \ .
\end{equation}
Together with the invariance of $K$, this implies in particular that
\begin{equation}\label{e.3-2a}
f(B_\eps(K)) \ \ssq \ B_\eps(K) \ .
\end{equation}
It follows that for any $(\theta,x)\in B_\eps(K)$
\begin{equation}
\label{e.3-2b}
  \|Df_\theta^n(x)\| \ \leq \ \alpha^n \quad \forall n\in\N \ .
\end{equation}
Consequently , we have
\begin{equation}
  \label{e.3-3a}
  x,x'\in K_\theta \textrm{ and } d(x,x')< 2\eps \quad \follows \quad
  d(f^n_\theta(x),f^n_\theta(x')) \ \leq \ \alpha^n \cdot d(x,x') \ \forall n\in\N \ .
\end{equation}

We now proceed in 4 steps. \medskip

\noindent {\em Step 1: \quad $K$ intersects every fibre in a finite number of points.}
\smallskip

Let $K_\theta:=\{x\in X:(\theta,x)\in K\}$.  As $K$ is compact, there exist
$(\theta_1,x_1) \ld (\theta_m,x_m)$ such that
\begin{equation}
  \label{e.3-4}
  K \ \ssq \ \bigcup_{k=1}^m B_\eps(\theta_k,x_k) \ .
\end{equation}
We will show that for any $\theta\in\Theta$ the cardinality of $K_{\theta}$,
denoted by $\#K_{\theta}$, is at most $m$.

Suppose for a contradiction that there exists $\theta_0\in\Theta$ with
$\#K_{\theta_0} > m$. We choose $m+1$ distinct points $\xi_1 \ld \xi_{m+1} \in
K_{theta_0}$ and let
\[
a \ = \ \min_{i\neq j} d(\xi_i,\xi_j) \ .
\]
Further, we fix $n\in\N$ such that $2\eps \cdot \alpha^{n}<a$ and choose, for
each for $i=1\ld m+1$, some $\xi_i'\in
\left(f^n_{\omega^{-n}(\theta_0)}\right)^{-1}\{\xi_i\} \in K$ (note that such
$\xi_i'$ exist since $f(K)=K$ and therefore
$f^n_{\omega^{-n}(\theta_0)}(K_{\omega^{-n}(\theta_0)})=K_{\theta_0}$). Due to (\ref{e.3-4}),
there exist $l\in\{1\ld m\}$ and $i,j\in\{1\ld m+1\}$ such that $\xi_i'$ and
$\xi_j'$ both belong to $B_\eps(x_l)$. Hence, the distance between the two
points is less than $2\eps$. Using (\ref{e.3-3a}) we conclude that
\begin{equation} \label{e.3-4a}
d(\xi_i,\xi_j) \ = \ d\left(f^n_{\omega^{-n}(\theta_0)}(\xi_i'),f^n_{\omega^{-n}(\theta_0)}(\xi_j')\right) \ \leq \
\alpha^n\cdot 2\eps \ < \ a \ ,
\end{equation}
contradicting the definition of $a$.
\medskip

\noindent {\em Step 2: \quad $\#K_\theta$ is constant on $\Theta$.}
\smallskip

We let
\[
k \ := \ \min_{\theta\in\Theta} \# K_\theta
\]
and fix $\theta_0$ with $\# K_{\theta_0} = k$. Suppose there exists $\theta \in
\Theta$ with $\# K_\theta > k$. Similar as in Step 1, we choose points $\xi_1
\ld \xi_{k+1} \in K_\theta$, let $a \ = \ \min_{i\neq j} d(\xi_i,\xi_j)$ and fix
$n_0\in\N$ such that $\alpha^n\cdot 2\eps < a \ \forall n\geq n_0$. Due to the
compacity of $K$, there exists $\delta > 0$ such that
\begin{equation} \label{e.3-5} K_{\theta'} \ \ssq \ B_\eps(K_{\theta_0})\quad
  \forall \theta' \in B_\delta(\theta_0) \ .
\end{equation}
By the minimality of $\omega$ on $\Theta$, there exists $n\geq n_0$
with $\omega^{-n}(\theta) \in B_\delta(\theta_0)$, such that
$K_{\omega^{-n}(\theta)} \ssq B_\eps(K_{\theta_0})$. However, as
$K_{\theta_0}$ only consists of $m$ points, at least two of the points
$\xi_1\ld \xi_{m+1}$, say $\xi_i$ and $\xi_j$, must have preimages
$\xi_i'$ and $\xi_j'$ under $f^n_{\omega^{-n}(\theta)}$ such that
$d(\xi_i',\xi_j') < 2\eps$.  Using (\ref{e.3-3a}) again we obtain
\begin{equation} \label{e.3-6} d(\xi_i,\xi_j) \ = \
  d\left(f^n_{\omega^{-n}(\theta)}(\xi_i'),f^n_{\omega^{-n}(\theta)}(\xi_j')\right) \ \leq
  \ \alpha^n\cdot 2\eps \ < \ a \ ,
\end{equation}
contradicting the definition of $a$.
\medskip

\noindent
{\em Step 3: \quad The distance between distinct points in $K_\theta$ is at least $2\eps$.}
\smallskip

The proof of this step is almost completely identical to that of Step 2. If
there exists $\theta_0\in\Theta$ such that two points in $K_{\theta_0}$ have
distance less than $2\eps$, then for any $n$ with $\omega^{-n}(\theta)$
sufficiently close to $\theta_0$ at least two of the $k$ points in $K_\theta$
will have preimages that are $2\eps$-close. Choosing $n$ sufficiently large and
using (\ref{e.3-3a}) once more, this leads to a contradiction in the same way as in
(\ref{e.3-4a}) and (\ref{e.3-6}).\medskip

\noindent
{\em Step 4: \quad The mapping $\theta \mapsto K_\theta$ is continuous.}
\smallskip

Fix $\theta_0\in\Theta$. We have to show that given any $\gamma>0$ there exists
$\delta>0$ such that $d(\theta,\theta_0) < \delta$ implies
$d_H(K_\theta,K_{\theta_0}) < \gamma$, where $d_H$ denotes the Hausdorff
distance on the space of subsets of $X$.

We may assume without loss of generality that $\gamma < \eps$. Due to the
compacity of $K$, there exists $\delta > 0$ such that $d(\theta,\theta_0)<
\delta$ implies $K_\theta \ssq B_\gamma(K_{\theta_0})$. However, since
$K_{\theta}$ and $K_{\theta_0}$ consist of exactly $k$ points which are at least
$2\eps$ apart, there must be exactly one point of $K_{\theta}$ in the
$\gamma$-neighbourhood of any point in $K_{\theta_0}$. Thus, we obtain
$d_H(K_\theta,K_{\theta_0}) < \gamma$ as required. \qed \medskip

\section{Saddle-node bifurcations: deterministic forcing}  \label{DeterministicForcing}

We come to the deterministic counterpart of Theorem~\ref{t.random-sn}.

\begin{theorem}[Saddle-node bifurcations, deterministic
  forcing] \label{t.deterministic-sn} Let $\omega$ be a homeomorphism of a
  compact metric space $\Theta$ and suppose that $(f_\beta)_{\beta\in [0,1]}$ is
  a parameter family of $\omega$-forced monotone ${\cal C}^2$-interval
  maps. Further, assume that there exist continuous functions $\gamma^-,\gamma^+
  : \Theta \to X$ with $\gamma^-<\gamma^+$ such that the following holds (for
  all $\beta\in[0,1]$ and $\theta\in\Theta$ where applicable).  \romanlist
\item[(d1)] There exist two distinct continuous $f_0$-invariant graphs and no
  $f_1$-invariant graph in $\Gamma$;
 \item[(d2)] $f_{\beta,\theta}(\gamma^\pm(\theta)) \geq
   \gamma^\pm(\omega(\theta))$;
 \item[(d3)] the maps $(\beta,\theta,x)\mapsto \partial_x^i f_\beta(\theta,x)$
   with $i=0,1,2$ and $(\beta,\theta,x)\mapsto \partial_\beta f_\beta(\theta,x)$
   are continuous;
 \item[(d4)] $f'_{\beta,\theta}(x) > 0$ for all $x\in \Gamma_\theta$;
 \item[(d5)] $\partial_\beta f_{\beta,\theta}(x) > 0 \ \forall x \in \Gamma_\theta$;
   \item[(d6)] $f''_{\beta,\theta}(x) > 0 \ \forall x\in \Gamma_\theta$;
 \listend

  Then there exists a unique critical parameter $\beta_c\in(0,1)$
  such that there holds:
\begin{itemize}
\item If $\beta < \beta_c$ then there exist two continuous $f_\beta$-invariant graphs
  $\varphi_\beta^-<\varphi^+_\beta$ in $\Gamma$. For any $\omega$-invariant measure $\mu$ we
  have $\lambda_\mu(\varphi^-_\beta)<0$ and $\lambda_\mu(\varphi^+_\beta)>0$.
\item If $\beta=\beta_c$ then either there exists exactly one continuous
  $f_\beta$-invariant graph $\varphi_\beta$ in $\Gamma$, or there exist two
  semi-continuous and weakly pinched $f_\beta$-invariant graphs $\varphi^-_\beta\leq
  \varphi^+_\beta$ in $\Gamma$, with $\varphi^-_\beta$ lower and
  $\varphi^+_\beta$ upper semi-continuous. If $\mu$ is an $\omega$-invariant
  measure then in the first case $\lambda_\mu(\varphi_\beta)=0$.  In the second
  case $\varphi^-_\beta(\theta)=\varphi^+_\beta(\theta)$ $\mu$-a.s.\ implies
  $\lambda_\mu(\varphi^\pm_\beta)=0$, whereas
  $\varphi^-_\beta(\theta)<\varphi^+_\beta(\theta)$
  $\mu$-a.s.\ implies $\lambda_\mu(\varphi^-_\beta)< 0$ and
  $\lambda_\mu(\varphi^+_\beta)>0$.
\item If $\beta > \beta_c$ then no $f_\beta$-invariant graphs exist in $\Gamma$.
\end{itemize}
\end{theorem}

\begin{remarks} \label{r.deterministic-forcing} \alphlist
\item In the above setting, we do not speak of equivalence classes of invariant
  graphs as in Section~\ref{RandomForcing}, but require invariant graphs to be
  defined everywhere. This results in a non-uniqueness of the invariant graphs
  in the above statement. For example, if $\omega$ has a wandering open set $U$,
  then the invariant graphs can easily be modified on the orbit of $U$. However,
  uniqueness can be achieved by requiring $\varphi_\beta^-$ to be the lowest
  and $\varphi^+_\beta$ to be the highest invariant graph in $\Gamma$.
\item Continuity and compacity imply that the derivatives in (d4)--(d6) are bounded
  away from zero by a uniform constant. In addition, if $\omega$ is minimal then it
  suffices to assume strict inequalities only for a single $\theta\in\Theta$, since
  for a suitable iterate the inequalities will be strict everywhere.
\item Again, a symmetric version holds for concave fibre maps (compare
  Remark~\ref{r.random-forcing}(d)).
\item We have to leave open here whether weakly pinched, but not pinched
  invariant graphs may occur at the bifurcation point in the above
  setting. While weakly pinched, but not pinched invariant graphs can be produced
  easily in general forced monotone maps, we conjecture that the additional
  concavity assumption excludes such behaviour in our setting.
\item The above result can be seen as a generalisation of results by
  the Alonso and Obaya \cite{alonso/obaya:2003} and Nunez and Obaya
  \cite{nunez/obaya:2007}, although the methods of proof are quite
  different. We discuss the relations in more detail in the next
  section.  \listend
\end{remarks}
\proof[Proof of Theorem~\ref{t.deterministic-sn}] As $f$ and $\gamma^\pm$ are
continuous, the sequences $\gamma^\pm_{\beta,n}$ defined by (\ref{e.gamma_n})
consist of continuous curves. Consequently, if the limits $\varphi^-_\beta$ and
$\varphi^+_\beta$ exist then due to the monotone convergence they are lower and
upper semi-continuous, respectively. Further, the sequences
$\gamma^\pm_{\beta,n}$ remain bounded in $\Gamma$ if and only if there exists an
$f_\beta$-invariant graph in $\Gamma$. In this case, $\varphi^-_\beta$ is the
lowest and $\varphi^+_\beta$ is the highest $f_\beta$-invariant graph in
$\Gamma$. We let
\begin{equation}
  \label{e.det-beta_c}
  \beta_c \ = \ \sup\left\{\beta\in[0,1]\left| \ \forall \beta'<\beta \ \exists \ 2
      \textrm{ uniformly separated } f_{\beta'}\textrm{-invariant graphs in } 
     \Gamma \right.\right\} \ .
\end{equation}
Note that we have $\beta_c \leq \beta_\mu$ for all $\omega$-invariant measures
$\mu$ (where $\beta_\mu$ is the critical parameter from
Theorem~\ref{t.random-sn}), since a pair of uniformly separated invariant graphs
is certainly  $\mu$-uniformly separated as well. \medskip

\noindent \underline{$\beta<\beta_c$:} \quad We have to show that
$\varphi^-_\beta$ and $\varphi^+_\beta$ are continuous, the statement about the
Lyapunov exponents then follows from Theorem~\ref{t.convexity}. As the two
graphs are uniformly separated, there exists $\delta>0$ such that
$\varphi^-_\beta(\theta)\leq \varphi^+_\beta(\theta)-\delta \ \forall
\theta\in\Theta$. Consequently, the point set $\Phi^-_\beta$ is contained in
$[\varphi^-_\beta,\varphi^+_\beta-\delta]$, and therefore the same is true for
the set $K:=\left[\varphi^-_\beta,\varphi^{-+}_\beta\right]$. Hence $K\cap
  \Phi^+_\beta =\emptyset$.

  Suppose $\mu$ is an $\omega$-invariant measure and $\varphi$ is an
  $(f_\beta,\mu)$-invariant graph contained in $K$. As there can be at most two
  $(f_\beta,\mu)$-invariant graphs in $\Gamma$ by Theorem~\ref{t.convexity}, we
  must have $\varphi=\varphi^-_\beta$ or $\varphi=\varphi^+_\beta$
  $\mu$-a.s.~. However, as $K\cap \Phi^+_\beta=\emptyset$ the case
  $\varphi=\varphi^+_\beta$ $\mu$-a.s.\ is not possible, such that
  $\varphi=\varphi^-_\beta$ $\mu$-a.s.~. Thus we have
  $\lambda_\mu(\varphi)=\lambda_\mu(\varphi^-_\beta)<0$ by
  Theorem~\ref{t.convexity}.

  Since $\mu$ and $\varphi$ were arbitrary, $K$ satisfies the assumptions of
  Lemma~\ref{l.curve-criterium} and we conclude that $K=\Phi^-_\beta$ is a
  continuous curve. Replacing $f$ with $f^{-1}$, which changes the signs of the
  Lyapunov exponents, the same argument shows that $\varphi^+_\beta$ is
  continuous as well. \medskip

  \noindent \underline{$\beta=\beta_c$ and $\beta>\beta_c$:} \quad Here the
  arguments are exactly the same as in the proof of Theorem~\ref{t.random-sn},
  with $(f,\mu)$-invariance replaced by $f$-invariance and measurable pinching
  by weak pinching.  \qed

  \medskip\medskip As in Section~\ref{RandomForcing}, we close with a discussion
  of bifurcations that take place on invariant subsets. If $M$ is a compact
  $\omega$-invariant subset of $\Theta$, then Theorem \ref{t.deterministic-sn}
  holds for the deterministic forcing system $(M,{\cal B},\omega_{|M})$
  and the parameter family $f_{\beta|M\times X}$ with new bifurcation
  parameter
\begin{equation*}
  \beta_c^M = \sup\left\{\beta\in[0,1] \left| \ \forall \beta'<\beta \
 \exists \textrm{ 2 uniformly separated } f_{\beta|M\times X}
 \textrm{-invariant graphs}\right.\right\}.
\end{equation*}
Obviously, we have
\begin{lemma} \label{l.subset-bif}
  Let $M\subseteq\Theta$ be compact and $\omega$-invariant. Then
  $\beta^M_c\geq\beta_c$.
\end{lemma}
 With the same notation as introduced
after Remark~\ref{r.parameter-monotonicity}, we have the following analogues to
Lemma~\ref{l.graph-monotonicity} and Lemma~\ref{l.bounded-orbits}.
\begin{lemma} \label{l.deterministic-graph-monotonicity}The function
  $\beta\mapsto\varphi_\beta^-(x)$ is increasing and the function
  $\beta\mapsto\varphi_\beta^+(x)$ is decreasing, for all
  $x\in\Gamma_{\theta}$, $\theta\in\Theta$.
\end{lemma}
We define $K(\beta)$ and $B(\beta)$ in the same way as in
Lemma~\ref{l.bounded-orbits}.
\begin{lemma} \label{l.gamma-bounded-orbits}
The following hold for all $\forall \beta\in[0,1]$.
\item[\, (i)] $K(\beta)$ is compact and $f_\beta$-invariant, $B(\beta)$ is
  compact and $\omega$-invariant.
\item[\, (ii)]
$K(\beta)=[\varphi_\beta^-,\varphi_\beta^+]$.
\item[\, (iii)]
If $\beta>\beta'$, then $B(\beta)\subseteq B(\beta')$ and
$K(\beta)\subseteq K(\beta')$
\end{lemma}
\proof The proof is identical to that of Lemma~\ref{l.bounded-orbits}, compacity
in (i) being a direct consequence of continuity.  \qed\medskip

As in Section~\ref{RandomForcing}, we can define a last bifurcation parameter
\[
\hat\beta_c \ = \ \sup\{\beta\in[0,1] \mid K(\beta) \neq \emptyset \} \
\]
and a bifurcation interval $I_c=[\beta_c,\hat\beta_c]$ over which the set of
$\Gamma$-bounded orbits vanishes. In contrast to the measurable setting, where
$K(\hat\beta_\mu)$ may be empty, we have
\begin{lemma}
  $K(\hat\beta_c) \neq \emptyset$.
\end{lemma}
\proof Due to Lemma~\ref{l.gamma-bounded-orbits}(iii) the sets
$K_n:=K(\hat\beta_c+1/n)$ form a nested sequence of compact sets. Hence
$K=\bigcap_{n\in\N} K_n$ is compact and non-empty, and continuity implies
$K=K(\hat\beta_c)$. \qed\medskip

In the minimal case, the bifurcation interval degenerates to a unique
bifurcation point.
\begin{lemma}
  If $\omega$ is minimal, then $B(\beta)=\Theta$ for $\beta\leq\beta_c$, and
  $B(\beta)=\emptyset$ for $\beta>\beta_c$.
\end{lemma}

Finally, we note that even if $\omega$ is uniquely ergodic with unique
invariant measure $\mu$, $\beta_c$ and $\beta_\mu$ need not coincide.
More precisely, we have $\beta_c \leq \beta_\mu$, but
$\beta_c<\beta_\mu$ may happen.

\section{Application to continuous-time systems}\label{ContinuousTime}

We now consider skew product flows
\[
\Xi_\beta \ : \ \R \times \Theta \times X \to \Theta\times X \quad , \quad (t,\theta,x) \mapsto
(\omega_t(\theta),\xi_\beta(t,\theta,x))
\]
generated by non-autonomous scalar differential equations
\[
x'(t) \ = \ F_\beta(\omega_t(\theta),x(t))
\]
with parameter $\beta\in[0,1]$ and base flow $\omega:\R\times\Theta\to\Theta$. We
concentrate on the deterministic case where $\Theta$ is a compact metric space
and $\omega:\R\times\Theta \to \Theta$ is a continuous flow. The random case can
be treated in a similar way.

Fix $t_0>0$ and let $f_\beta(\theta,x):=\Xi_\beta(t_0,\theta,x)$. We say
$\varphi:\Theta\to X$ is a {\em $\Xi_\beta$-invariant graph} if
$\xi_\beta(t,\theta,\varphi(\theta))=\varphi(\omega_t(\theta)) \ \forall t
\in\R ,\theta\in \Theta$. Obviously, in this case $\varphi$ is a
$f_\beta$-invariant graph as well.  Let $\gamma^-,\gamma^+ : \Theta\to X$ be
${\cal C}^1$-functions and suppose that \romanlist
\item[$(c1)$] there exist two $\Xi_0$-invariant graphs but no
  $\Xi_1$-invariant graph in $\Gamma$;
\item[$(c2)$] $\partial_t \gamma^\pm(\omega_t(\theta)) \ \leq \
  F_\beta(\omega_t(\theta),\gamma^\pm(\omega_t(\theta))) \ \forall t\in\R,
  \ \theta\in\Theta$ and $\beta\in[0,1]$; \listend We will see
  below that in the situation we consider this implies assumption $(d1)$ from
  Theorem~\ref{t.deterministic-sn} for $f_\beta$.  Moreover, due to $(c2)$ the
  map $t\mapsto \xi_\beta(t,\theta,\gamma^\pm(\theta)) -
  \gamma^\pm(\omega_t(\theta))$ is either strictly positive or zero and
  non-decreasing, and therefore non-negative for all $t>0$.  Consequently
\begin{equation}\label{e.flow-boundaries}
\xi_\beta(t,\theta,\gamma^\pm(\theta)) \ \geq \
\gamma^\pm(\omega_t(\theta)) \quad \forall t\in\R^+,\ \theta\in\Theta \ .
\end{equation}
Further, assume that
\romanlist
\item[$(c3)$] $(\beta,\theta,x) \mapsto F_\beta(\theta,x)$,
  $(\beta,\theta,x)\mapsto \partial_xF_\beta(\theta,x)$ and
  $(\beta,\theta,x)\mapsto \partial_\beta F_\beta(\theta,x)$ are continuous;
  \listend Then
  $\partial_xf_{\beta,\theta}(x),\ \partial_x^2f_{\beta,\theta}(x)$ and
  $\partial_\beta f_{\beta,\theta}(x)$ exist and are continuous. More
  explicitly, we have the following formulae.
  \begin{eqnarray}\label{e.F1}
    \partial_x f_{\beta,\theta}(x) & = & \exp\left(\int_0^{t_0} \partial_x
  F_\beta(\omega_s(\theta),\xi_\beta(s,\theta,x))\ ds \right) \\ \label{e.F2}
\partial_x^2 f_{\beta,\theta}(x) & = & \exp\left(\int_0^{t_0} \partial_x
  F_\beta(\omega_s(\theta),\xi_\beta(s,\theta,x))\ ds \right) \cdot\\ & &\nonumber
\int_0^{t_0} \partial_x^2 F_\beta(\omega_s(\theta),\xi_\beta(s,\theta,x))
\cdot \partial_x \xi_\beta(s,\theta,x) \ ds \ .\\\label{e.F3}
 \partial_\beta f_{\beta,\theta}(x)  & = & \int_0^{t_0} \partial_\beta
 F_\beta(\omega_s(\theta),\xi_\beta(s,\theta,x)) \cdot \int_s^{t_0}\partial_x
 F_\beta(\omega_r(\theta),\xi_\beta(r,\theta,x)) \ dr \ ds \ .
  \end{eqnarray}
From (\ref{e.F1}), we see that
\romanlist
\item[$(c4)$] $\partial_x F_\beta(\theta,x) > 0 \ \forall (\theta,x,\beta)
  \in\Theta\times X \times [0,1]$ \listend implies $\partial_x f_{\beta,\theta}
  > 0$ and hence (d4).
From (\ref{e.F2}) we can deduce that
\romanlist
\item[$(c5)$] $\partial_\beta F_\beta(\theta,x) > 0 \ \forall (\theta,x,\beta)
  \in\Theta\times X \times [0,1]$ \listend implies $\partial_\beta
  f_{\beta,\theta}(x) >0$, such that (d5) holds.
 Finally \romanlist
\item[$(c6)$] $\partial_x^2 F_\beta(\theta,x) > 0 \ \forall (\theta,x,\beta)
  \in\Theta\times X \times [0,1]$ \listend yields the strict convexity of
  $f_{\beta,\theta}$, such that (d6) holds. \medskip

  Now suppose, that for some $\beta\in[0,1]$ the flow $\Xi_\beta$ has two
  invariant graphs in $\Gamma$. These can be obtained as the monotone limits of
  the sequences
\[
\gamma^-_{\beta,t}(\theta) =
\xi_\beta(t,\omega_{-t}(\theta),\gamma^-(\omega_{-t}(\theta))) \qquad
\textrm{and} \qquad \gamma^+_{\beta,t}(\theta) =
\xi_\beta(-t,\omega_t(\theta),\gamma^+(\omega_t(\theta))) \ ,
\]
by taking
\[
\varphi^-_\beta(\theta) = \tLim \gamma^-_{\beta,t}(\theta) \qquad \textrm{and} \qquad
\varphi^+_\beta(\theta) = \tLim \gamma^+_{\beta,t}(\theta) \ .
\]
Since these are also $f_\beta$-invariant, $f_0$ has two invariant
graphs in $\Gamma$.

Conversely, if $f_\beta$ has an invariant graph $\varphi$ in $\Gamma$, then for
all $\theta\in\Theta$ and $t\in\R$ the points
$\Xi_\beta(t,\theta,\varphi(\theta))$ remain in $\Gamma$. (Note that due to
the monotonicity of the flow in the fibres and (\ref{e.flow-boundaries}), orbits
which have left $\Gamma$ can never return.) Hence, the graphs of
$\gamma^\pm_{\beta,t}$ remain in $\Gamma$ for all $t$ and therefore
$\Xi_\beta$ has invariant graphs $\varphi^-_\beta$ and $\varphi^+_\beta$ as
well (which might coincide). Consequently, if $\Xi_\beta$ has no invariant
graphs, then the same is true for $f_\beta$. This shows that (c1) implies
(d1) and altogether that (c1)--(c6) imply (d1)--(d6). This leads to the
following continuous-time version of Theorem~\ref{t.deterministic-sn}, which is
a generalisation of results in \cite{novo/obaya/sanz:2004,nunez/obaya:2007} on
strictly ergodically forced convex scalar differential equations.
\begin{thm}
  Suppose $(F_\beta)_{\beta\in[0,1]}$ satisfies (c1)--(c6). Then there exists a
  unique critical parameter $\beta_c\in(0,1)$, such that
\begin{itemize}
\item If $\beta < \beta_c$ then there exist two continuous $\Xi_\beta$-invariant graphs
  $\varphi_\beta^-<\varphi^+_\beta$ in $\Gamma$. For any $\omega$-invariant measure $\mu$ we
  have $\lambda_\mu(\varphi^-_\beta)<0$ and $\lambda_\mu(\varphi^+_\beta)>0$.
\item If $\beta=\beta_c$ then either there exists exactly one continuous
  $\Xi_\beta$-invariant graph $\varphi_\beta$ in $\Gamma$, or there exist two
  semi-continuous and weakly pinched $\Xi_\beta$-invariant graphs $\varphi^-_\beta\leq
  \varphi^+_\beta$ in $\Gamma$, with $\varphi^-_\beta$ lower and
  $\varphi^+_\beta$ upper semi-continuous. If $\mu$ is an $\omega$-invariant
  measure then in the first case $\lambda_\mu(\varphi_\beta)=0$.  In the second
  case $\varphi^-_\beta(\theta)=\varphi^+_\beta(\theta)$ $\mu$-a.s.\ implies
  $\lambda_\mu(\varphi^\pm_\beta)=0$, whereas
  $\varphi^-_\beta(\theta)<\varphi^+_\beta(\theta)$
  $\mu$-a.s.\ implies $\lambda_\mu(\varphi^-_\beta)< 0$ and
  $\lambda_\mu(\varphi^+_\beta)>0$ otherwise.
\item If $\beta > \beta_c$ there exist no $\Xi_\beta$-invariant graphs in
  $\Gamma$.
\end{itemize}

\end{thm}

\section{Some examples} \label{Examples}

In this section, the preceding results in this article will be illustrated by some
explicit examples.  In order to start with a simple case, we first
choose the base transformation $\omega$ to be an irrational rotation of
the circle, that is, $\omega: \kreis\to\kreis,\ \theta\mapsto\theta+\rho
\bmod 1$, where $\rho$ is the golden mean. Then minimality of $\omega$
and ergodicity of the Lebesgue measure $\mu$ on $\kreis$ will imply
that the bifurcation parameters $\beta_{\mu}$ and $\beta_c$ for the
measure-theoretic and the topological setting coincide, and that no
additional bifurcation parameters in the sense of
Remark~\ref{r.parameter-monotonicity} and Lemma~\ref{l.subset-bif}
exist. Further, it is well-known that a suitable choice of the fibre
maps $f_{\beta,\theta}$ will lead to a non-smooth bifurcation, in the
sense that a pair of non-continuous pinched invariant graphs exists at
the bifurcation point (instead of a single neutral and continuous
curve). In this context, these graphs are usually called strange
non-chaotic attractors, respectively repellers, depending on the sign
of the Lyapunov exponent
\cite{grebogi/ott/pelikan/yorke:1984,jaeger:2006a}.
 
In order to obtain such a non-smooth bifurcation, we choose
\begin{equation} \label{e.first-example}
f_{\beta}(\theta,x) \ = \ (\omega(\theta),\arctan(\alpha x) - 2\beta -
g(\theta)) \ ,
\end{equation}
where $g(\theta)=(\sin(2\pi \theta)+1)/2$. In fact, in order to apply
rigorous results on the existence of strange non-chaotic attractors a
slightly different choice of the forcing function would be required,
since such results are still due to a number of technical constraints
\cite{jaeger:2006a}. However, for the pictures obtained by simulations
there is hardly any difference. For the application of our results to
this parametrised family, we will use one of the analogue versions of
Theorem \ref{t.random-sn}, respectively
Theorem~\ref{t.deterministic-sn}, mentioned in
Remarks~\ref{r.random-forcing}(d) and
\ref{r.deterministic-forcing}(c). More precisely, instead of convexity
in $(r7)$ and $(d6)$ we will require concavity and instead of positive
derivative with respect to $\beta$ in $(r6)$ and $(d5)$ we will
require negative derivative. In $(r2)$ and $(d2)$ the inequalities
then need to be reversed. All other conditions remain as before, and
the only difference in the statement is that the signs of the Lyapunov
exponents will be reversed.

\begin{figure}[h!]
\begin{center}
\includegraphics[scale=1.19]{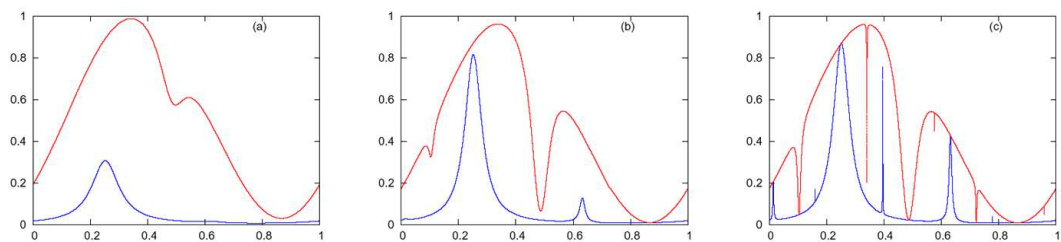}
\caption{\label{fig1} \small Invariant graphs for the $1$-parameter family
  $f_\beta(x,y)=(x+\rho, \arctan(\alpha y)-2\beta-\gamma(\sin(2\pi
  x)+1))$ where $\rho$ is the golden mean, $\alpha=100$, $\gamma=1/2$,
  and (a) $\beta=0.265$, (b) $\beta=0.275$, (c) $\beta=0.2753743$.}
\end{center} 
\end{figure}

For all $\beta\geq 0$, the curves $\gamma^-\equiv 0$ and
$\gamma^+\equiv 2$ satisfy $f^{\pm 1}_{\beta*}\gamma^\pm \leq
\gamma^\pm$.  Conditions $(r3)$--$(r7)$ and $(d3)$--$(d6)$ are
obviously verified. In order to check $(r1)$, respectively $(d1)$,
note that for all sufficiently large $\alpha$ (say, $\alpha\geq 20$),
the curve $\psi$ given by
$\psi(\theta)=\dreiviertel-\halb\sin(2\pi(\theta-\rho))$ satisfies
$f_{0*}\psi \geq \psi$. As argued in the proof of
Theorem~\ref{t.random-sn}, this implies the existence of two
$f_0$-invariant graphs (compare (\ref{e.intermediate-graph})), whereas
the non-existence of $f_{\beta_1}$-invariant graphs in $\Gamma$ can be
seen from the fact that $f_{1,0}(2)<0$. Consequently
(\ref{e.first-example}) satisfies all assumptions of (the analogue
version of) Theorems~\ref{t.random-sn} and \ref{t.deterministic-sn},
and we obtain the existence of a saddle-node bifurcation in $\Gamma$.
Figure~\ref{fig1} shows the approach of the upper and lower invariant
graph in $\Gamma$. In (c), $\beta=0.2753743$ is a good approximation
of the bifurcation point and the picture gives an idea of the strange
non-chaotic attractor-repeller pair that emerges.

For slightly larger parameters $\beta$, the invariant graphs in
$\Gamma$ disappear. In this case, all trajectories converge to an
attracting continuous invariant graph, in the region below
$\kreis\times\{0\}$, which exists throughout the whole parameter
range. \medskip

In order to construct an example with a more complex bifurcation
pattern, in the sense discussed at the end of Sections
\ref{RandomForcing} and \ref{DeterministicForcing}, we need a base
transformation that exhibits more complicated dynamics and, in
particular, a multitude of invariant measures and minimal sets.
Evidently, the canonical choice is to use a two-dimensional
transformation, since this allows at the same time for the required
complex behaviour and a graphical representation of the invariant
graphs of the resulting three-dimensional system. Our choice is the
map
\begin{equation}\label{e.secondbase} \textstyle
\omega(\theta_1,\theta_2) \ = \ \left(\theta_1+
  \halb\sin\left(2\pi\left(\theta_2+\halb\sin(2\pi\theta_1)\right)\right),
  \theta_2+\halb\sin(2\pi\theta_1)\right)
\ ,
\end{equation}
which has been studied in its own right in the context of quantum
dynamics \cite{leboeuf:1990,geiseletal:1991}.

\begin{figure}[h!]
\begin{center} 
  \includegraphics[scale=1.22]{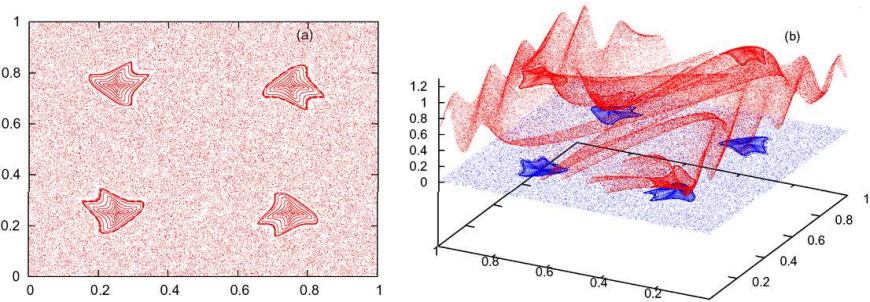}\caption{\label{fig2}\small (a) Phase
    portrait of the map $\omega$ from (\ref{e.secondbase}). (b) The
    two invariant graphs at the bifurcation point $\beta_c\simeq 0.1855650809$ for the
    parametrised family
    $f_\beta(\theta,x)=(\omega(\theta),f_{\beta,\theta}(x))$ with
    $\omega$ from (\ref{e.secondbase}) and $f_{\beta,\theta}$ defined
    by (\ref{e.secondexample}).}
\end{center}
\end{figure}

It is known that $\omega$ has both an uncountable number of invariant
ergodic measures and of minimal sets (this is due to the fact that
its rotation set has non-empty interior, see \cite{jaeger:2011} for a
discussion). For the illustration, it is particularly convenient that
$\omega$ exhibits four (star-shaped) elliptic islands, centred around
the points of two period-2 orbits
$M_1=\left\{\left(\viertel,\viertel\right),\left(\dreiviertel,\dreiviertel\right)\right\}$
and
$M_2=\left\{\left(\viertel,\dreiviertel\right),\left(\dreiviertel,\viertel\right)\right\}$ (see Figure~\ref{fig2}(a)).

As fibre maps, we choose
\begin{equation} \label{e.secondexample} \textstyle
  f_{\beta,\theta}(x) \ = \ \arctan(\alpha x)
  -2\beta-\gamma(\sin(2\pi\theta_1)\sin(2\pi\theta_2)+1) \ .
\end{equation}
Note that for $\gamma>0$ the $\theta$-dependent term
$-\gamma\sin(2\pi\theta_1)\sin(2\pi\theta_2)$ takes its global minimum
exactly at the two points of the two-periodic orbit $M_1$. This
implies that $M_1$ is the minimal set on which the first bifurcation
occurs, that is, $\beta_c^{M_1}=\beta_c<\beta_c^M \ \forall
\textrm{minimal sets } M\neq M_1$. Equivalently, $M_1$ is exactly the
set of points on which the two invariant graphs touch at the
bifurcation point. Furthermore, since
$f_{\beta,\left(\viertel,\viertel\right)} =
f_{\beta,\left(\dreiviertel,\dreiviertel\right)}$, the bifurcation
pattern of $f_{\beta|M_1}$ is the same as the one of the one-dimensional family
\[
g_\beta(x) \ = \ f_{\beta,\left(\viertel,\viertel\right)}(x) \ = \ \arctan(\alpha x) - 2\beta - 2\gamma \ .
\]
This allows to determine the precise bifurcation point, namely 
\begin{equation} \label{e.bif-point}
\beta_c \ = \ \halb\arctan(\sqrt{\alpha-1})-\frac{\sqrt{\alpha-1}}{2\alpha}-\gamma \ .
\end{equation}
For $a=100$ and $\gamma=1/2$ we obtain $\beta_c\simeq 0.1855650809$.

\begin{figure}[h!]
\begin{center}
\includegraphics[scale=1.3]{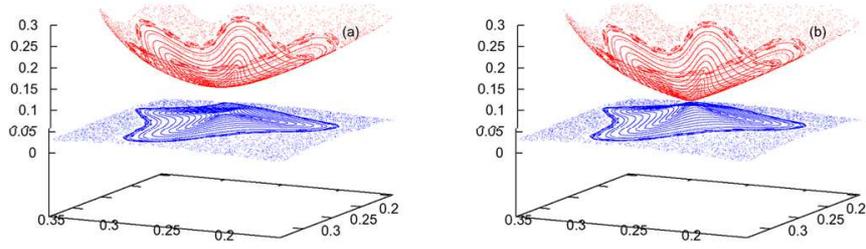}
\caption{\label{fig3} \small Closer view of the two invariant graphs
  over the islands `centred' at the period $2$ point $(1/4,1/4)$.
  (a) $\beta=\beta_c-0.005$, (b) $\beta=\beta_c$.}
\end{center}
\end{figure}

Figure~\ref{fig2}(b) shows the two invariant graphs in $\Gamma=\torus
\times [0,2]$ at this bifurcation point. The validity of the
assumptions of Theorems \ref{t.random-sn} and
\ref{t.deterministic-sn} is checked in a similar way as in the
previous example. The picture becomes clearer in Figure~\ref{fig3}\,
where the restriction of the two invariant graphs over a neighbourhood
of $\left(\viertel,\viertel\right)$ is plotted, slightly before the
bifurcation point in (a) and at the bifurcation point in (b).

\begin{figure}[h!]
\begin{center}
\begin{tabular}{cc}
\includegraphics[scale=1.27]{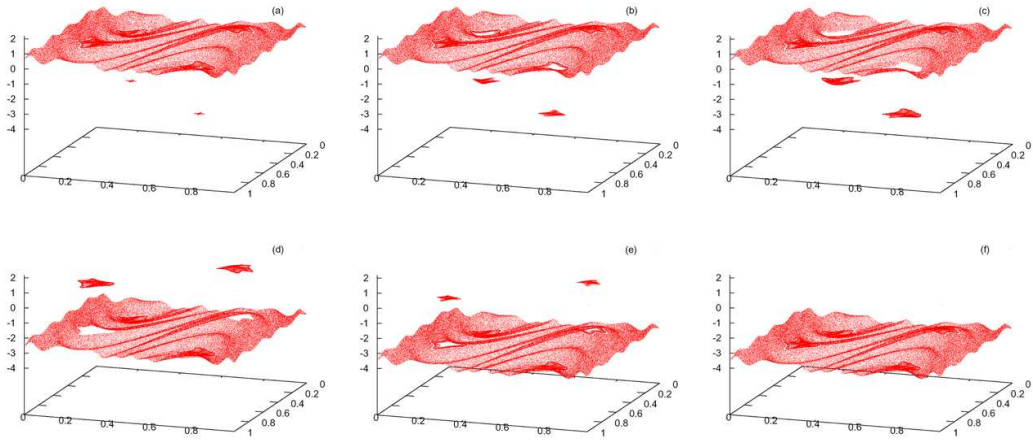}
\end{tabular}\caption{\label{fig4}\small Invariant graphs for $f_\beta$ where $\alpha=100$,
  $\gamma=1/2$, and (a) $\beta=\beta_c+0.0005$, (b) $\beta=\beta_c+0.01$, 
  (c) $\beta=\beta_c+0.0269$, (d) $\beta=\beta_c+0.02725$, (e) $\beta=\beta_c+0.485$, 
  (f) $\beta=\beta_c+0.5$.}
\end{center}
\end{figure}

Similarly to the previous example, there exists a third invariant graph
below $\torus \times \{0\}$, which is continuous and attracting and
persists throughout the whole parameter range. Once the bifurcation
has taken place over a minimal set $M$, this graph attracts all
trajectories in $M\times [-5,2]$. Consequently, the upper bounding
graph $\varphi^+_M$ {\em `drops down'} from above 0 to below at the
bifurcation point $\beta_c^M$. This happens first for $M_1$, and
subsequently for all the invariant circles in the elliptic island,
starting in the middle and moving outwards (see
Figure~\ref{fig4}(a)--(c)).  Note that in all pictures in Figure~\ref{fig4}
only the upper bounding graph is plotted, for the sake of better
visibility.

When the outer boundary of the two elliptic islands containing $M_1$
is reached, the complement of the elliptic islands (the chaotic region
in the sense of \cite{jaeger:2011}) drops in one go. Finally, the
invariant circles over the remaining two elliptic islands drop down
one by one, in reversed order, moving inwards from the outside (note
that on $M_2$ the $\theta$-dependent term takes its global maximum).

\begin{figure}[h!]
\begin{center}
  \includegraphics[scale=1.1]{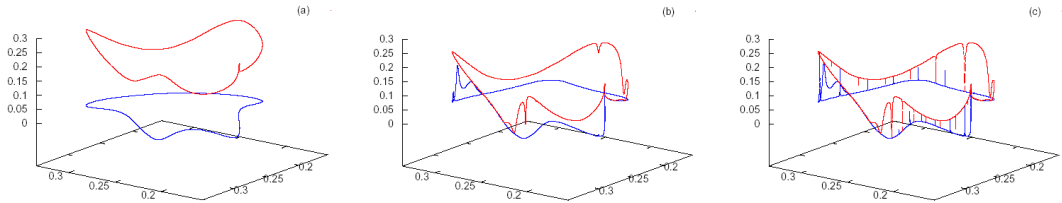}\caption{\label{fig5}\small
    Closer view of two invariant circles above the island centred at
    the point $(1/4,1/4)$. Here, $\alpha=200$, $\gamma=1$, and (a)
    $\beta=\beta_c+0.0035$, (b) $\beta=\beta_c+0.03516$, (c)
    $\beta=\beta_c+0.035164103$. $\beta_c$ is again determined by
    \eqref{e.bif-point}. Note that $\beta_c$ is negative in this case.
    Hence, strictly speaking a reparametrisation would be necessary to meet the formal
    requirements of Theorem~\ref{t.deterministic-sn}, but we omit the
    details.  }
\end{center}
\end{figure}

Finally, in Figure~\ref{fig5}, the bifurcation over one of the invariant
circles of the elliptic island is shown. Although embedded in
dimension two, the underlying dynamics are just those of an irrational
rotation. Consequently, from a qualitative point of view, the situation
is exactly the same as in the first example. Again, the
non-uniform approach of the invariant circles can be observed, which
is typical for the creation of strange non-chaotic attractors and
repellers at the bifurcation point.

\end{document}